\documentclass[journal]{IEEEtran}
\usepackage[short]{optidef}
\usepackage{graphicx, float}
\usepackage{amsmath}
\usepackage{siunitx}
\usepackage{amssymb}
\usepackage{amsthm}
\usepackage{booktabs}
\usepackage{amsfonts}
\usepackage{arydshln}
\usepackage{mathtools}
\usepackage{cite}
\usepackage[subtle]{savetrees}
\usepackage{soul}
\usepackage{color} 
\usepackage{flushend}
\usepackage[linesnumbered, ruled]{algorithm2e}
\definecolor{mypink1}{rgb}{0.858, 0.188, 0.478}

\newtheorem{theorem}{Theorem}[section]

\newtheorem*{remark}{Remark}

\newtheorem{definition}{Definition}[section]
\ifCLASSOPTIONcompsoc
    \usepackage[caption=true, font=normalsize, labelfont=sf, textfont=sf]{subfig}
\else
\usepackage[caption=false, font=footnotesize]{subfig}
\fi
\IEEEoverridecommandlockouts  

\title{Grid-aware aggregation and realtime disaggregation of distributed energy resources in radial networks
}

\author{Nawaf~Nazir$^\ast$, \IEEEmembership{Member, IEEE} and Mads~Almassalkhi$^\dagger$, \IEEEmembership{Senior Member, IEEE}
\thanks{
$^{\ast}$N. Nazir (nawaf.nazir@pnnl.gov) is affiliated with Pacific Northwest National Lab, Richland, WA 99354.
$^{\dagger}$M. Almassalkhi (malmassa@uvm.edu) is with the Department of Electrical and Biomedical Engineering, The University of Vermont, Burlington, VT 05405, USA. Support from the U.S. Department of Energy’s Advanced Research Projects Agency—Energy (ARPA-E) Award DE-AR0000694 is gratefully acknowledged. M. Almassalkhi was supported in part by the National Science Foundation (NSF) Award ECCS-2047306.}
}

\begin{document}
\maketitle
\begin{abstract}\label{abstract}
Dispatching a fleet of distributed energy resources (DERs) in response to wholesale energy market or regional grid signals  requires solving a challenging disaggregation problem when the DERs are coupled within a distribution network. This manuscript presents a computationally tractable convex inner approximation for the optimal power flow (OPF) problem that quantifies a feeder's aggregate DERs hosting capacity and enables a realtime, grid-aware control policy for DERs in radial distribution networks. The inner approximation is derived by considering convex envelopes on the nonlinear terms in the AC power flow equations. 
The resulting convex formulation is then used to derive provable nodal injection limits, such that any combination of DER dispatches within their respective nodal limits is guaranteed to be AC admissible. These nodal injection limits are then used to construct a realtime, open-loop control policy for dispatching DERs at each location in the network to deliver fast grid services in aggregate. The IEEE-37 distribution network is used to validate the technical results and illustrate use cases.
\end{abstract}

\section{Introduction}\label{sec:introduction}
The distribution system was engineered under the assumption that residential and commercial customers would only have power directed to them from the bulk grid. However, the increasing penetration of solar PV in distribution feeders has created so-called ``prosumers'' who (at times) can supply the grid with energy rather than just consume it. This leads to \textit{reverse power flows} that can result in unexpected violations of voltage and transformer limits, which may negatively impact system reliability~\cite{driesen2006distributed}. Furthermore, the significant variability inherent to solar PV generation  challenges traditional distribution system operating paradigms. In addition, with ubiquitous connectivity, behind-the-meter (BTM) smart appliances and DERs will soon underpin a demand that is inherently flexible.
Many works in literature, such as~\cite{hao2014aggregate,tindemans2015decentralized}, provide methodologies for aggregating the flexibility of DERs to provide grid services. The authors in~\cite{muller2017aggregation} employ transactive energy principles as way to disaggregate flexibility amongst the individual DERs. However, none of these methods consider the underlying network, which may become overloaded when flexible demand is deployed at scale in a distribution feeder.
The optimal power flow (OPF) represents a method for algorithms to improve reliability of the grid and responsiveness of flexible resources (e.g., batteries, PV inverters). However, due to the sub-minutely timescale of the solar PV variability, these algorithms must be computationally tractable and, yet, representative of the physics. That is, grid optimization algorithms should ensure admissible network operations~\cite{kristov2016tale}. 

Since Carpentier's original OPF formulation~\cite{carpentier1962contribution} and subsequent improvements in optimization solvers, the OPF problem has become a powerful methodology for optimizing the dispatch of various grid resources. This is because OPF-based methods can account for the underlying grid physics, static network constraints on voltages and apparent branch flows, and resource limitations. However, it was also recognized early on that the nonlinear AC power flow equations that model the underlying grid physics render the AC OPF problem non-convex\cite{molzahn2017computing}.
\begin{figure}[t]
\centering
\includegraphics[width=0.51\textwidth]{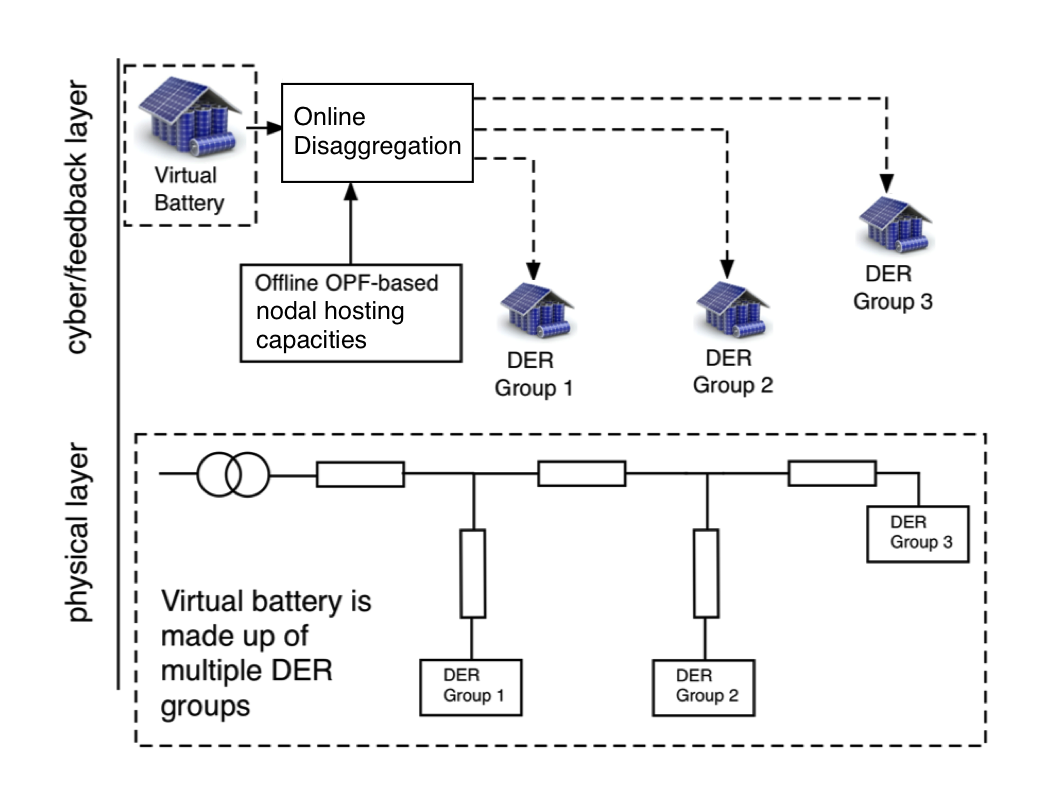}
\caption{\label{fig:cyber_physical_model} A schematic representation of the system model. The physical layer represents the circuit that couples the different DER groups into an aggregate virtual battery, whereas the cyber layer enables the disaggregation of the virtual battery market signal to the DER groups based on the nodal hosting capacities that are determined offline. VB image source: https://esdnews.com.au/}
\end{figure}
To overcome the computational challenges associated with non-convex OPF formulations, many recent techniques involve using linear approximations or convex relaxations \cite{molzahn_hiskens-fnt2019}. Traditional optimization techniques for dispatching resources include linear OPF-based \textit{LinDist} models~\cite{baran1989optimal}. These models work well when the distribution system is not heavily loaded (e.g., low losses). In~\cite{nazir2019convex}, it is shown how optimal DER dispatch based of the \textit{LinDist} model can lead to voltage violations. Similarly, the authors in~\cite{dvijotham2016error} quantified the errors associated with more general linear power flow approximations. Recently, improved linear approximations of the power flow equations have been proposed that provide improved accuracy over a wider range of operation~\cite{bolognani2015fast,arnold2016optimal}. However, the solution space of the AC power flow equations is highly non-convex, which means that linear approximations alone cannot guarantee a network-admissible dispatch under all (net) load conditions. 

Beyond linear approximations, recent attention in literature has focused on convex relaxations of the AC power flow equations, including second-order cone programs, semi-definite programs, and quadratic relaxations~\cite{taylor2015convex}.
Several works in literature have shown that, under certain analytical conditions, these relaxations can be exact and the solution of the relaxed convex problem then represents the global optimum of the original non-convex AC OPF problem~\cite{gan2015exact}. However, these conditions fail to hold under extreme solar PV injections, when the network experiences reverse power flows, which engenders a non-zero duality gap solution that may not be network admissible, i.e., not feasible in the original AC OPF formulation~\cite{huang2017sufficient}.

However, in many practical applications, guaranteeing network admissibility is often more valuable than finding the globally optimal solution. 
The authors in~\cite{molzahn2019grid}, develop an optimization-based method to certify whether a specific DER dispatch will result in constraint violations. However, they do not discuss whether a range of DER control actions is network-admissible. In~\cite{lee2019convex}, the authors employ a convex restriction that guarantees an AC admissible solution, which they utilize in~\cite{lee2020feasible} to determine an admissible path from a known initial operating point to a desired final operating point. {\color{black}However, their method requires careful coordination of generator outputs due to the ellipsoidal nature of their convex restriction. Note that this manuscript obtains a hyper-rectangle for the convex inner approximation, which enables DERs to operate independently.} 
In~\cite{nick2017exact}, the convex OPF formulation is based on an augmented second-order cone relaxation. 
The authors in~\cite{ROSS2020_PSCC} solve a large number of non-convex OPF problems to determine nodal injection bounds. However, these methods either rely on non-convex techniques or they cannot ensure that the full range of DER dispatch is network-admissible.

The manuscript herein presents a novel convex approximation of the AC OPF problem to quantify the network-admissible range of DER nodal injections in radial, balanced distribution feeders. In general, obtaining a convex inner approximation is NP-hard~\cite{chang1986polynomial}. However, the work herein uses the nonlinear branch-flow model (BFM) formulation of the AC power flow equations 
to define a convex envelope on the nonlinear terms relating the branch currents, nodal voltages, and apparent power flows. This is combined with the remaining linear elements of the BFM to form a convex inner approximation and ensures that all feasible (and, hence, optimal) solutions in the convex OPF are also feasible in the non-convex AC OPF formulation. We denote such a solution as \textit{network admissible} or \textit{AC admissible}.
From this approach, we achieve an OPF formulation that exhibits computational solve times similar to that of linear formulations with the added (and crucial benefit) that the formulation guarantees admissible solutions. This convex inner approximation is then utilized to determine the admissible nodal DER dispatch ranges for a radial network, i.e., any combination of dispatching DERs is network-admissible as long as each node is dispatched within its computed DER range. This methodology represents a major contribution in grid-aware dispatch of networked grid assets in distribution feeders and overcomes practical limitations of methods that rely on repeatedly solving centralized AC OPFs at each instant in receding-horizon fashion~\cite{nazir_inner} or require extensive, realtime communication of network and DER data~\cite{emiliano_VPP}. \textcolor{black}{The convex inner approximation methodology was first introduced in~\cite{nazir2019voltage}, but was utilized to optimize the reactive power set-points of controllable, discrete mechanical assets to maximize voltage margins (i.e., find a single optimal operating point). In this work, we have significantly extended the CIA approach to determine a range of admissible DER dispatch in radial feeders.} The main contributions of this manuscript are as follows:
\begin{itemize}
    \item Generalizes a CIA of the AC OPF problem that is applicable to any radial, balanced distribution feeder, such as those with a mix of inductive and capacitive branches and with branch current and nodal voltage limits. This improves over~\cite{nazir2019voltage}.
    \item The resulting CIA is employed to optimize the feeder's DER nodal hosting capacities, which represent the ranges of admissible injections for DERs at each node in the network such that all branch flows and nodal voltages are within limits (i.e., network admissible). Thus, the optimized DER nodal capacities can then be trivially aggregated to form the network's capacity for flexibility. Provable guarantees are provided for admissibility over the entire range of the nodal DER dispatch.
    \item 
    The CIA formulation provably guarantees existence and uniqueness of the underlying AC power flow solutions for the entire range of DER dispatches. This is achieved by adding another set of convex constraints to the CIA formulation.
    \item The admissible DER nodal capacities beget an open-loop, realtime disaggregation policy that accounts for network reliability constraints while enabling fast grid services.
    \item Simulation-based analysis assesses practicality of the proposed methods and investigates the conservativeness of the results and different reactive power strategies for augmenting a feeder's aggregate DER nodal capacity.

\end{itemize}

The remainder of the manuscript is organized as follows. Section~\ref{sec:Opt_form} develops the mathematical formulation of the convex inner approximation for the OPF problem using the robust bounds on nonlinear terms. Section~\ref{sec:Feas_guar} provides admissibility guarantees for the obtained DER nodal capacity and proposes an iterative algorithm that augments the admissible range. In Section~\ref{sec:Reac_power}, we present and analyze the effect of different nodal reactive power control policies to further augment the range of AC admissible flexibility that can be dispatched reiliably, whereas Section~\ref{sec:real_time} describes a realtime dispatch policy to disaggregate flexibility over a network in a grid-aware manner using the DER nodal capacities. Finally, Section~\ref{sec:conclusion} concludes the manuscript and lays out future research directions. 

\section{Formulating the Convex Inner Approximation}\label{sec:Opt_form}
The nonlinear \textit{DistFlow} model accurately represents the underlying physics for a radial, balanced AC distribution network~\cite{baran1989optimal}. However,  using \textit{DistFlow} in an AC OPF setting results in a non-convex optimization problem. Common techniques that employ linear approximations or convex relaxations are only valid under certain technical assumptions or near a pre-defined operating point.
In this section, we develop a novel convex inner approximation of the AC OPF that is used to compute the range of allowable nodal net injections, such that any combination of nodal injections within those ranges are guaranteed to satisfy AC limits for voltages and branch flows.

\subsection{Mathematical model}\label{sec:math_model}
\begin{figure}[t]
\centering
\includegraphics[width=0.38\textwidth]{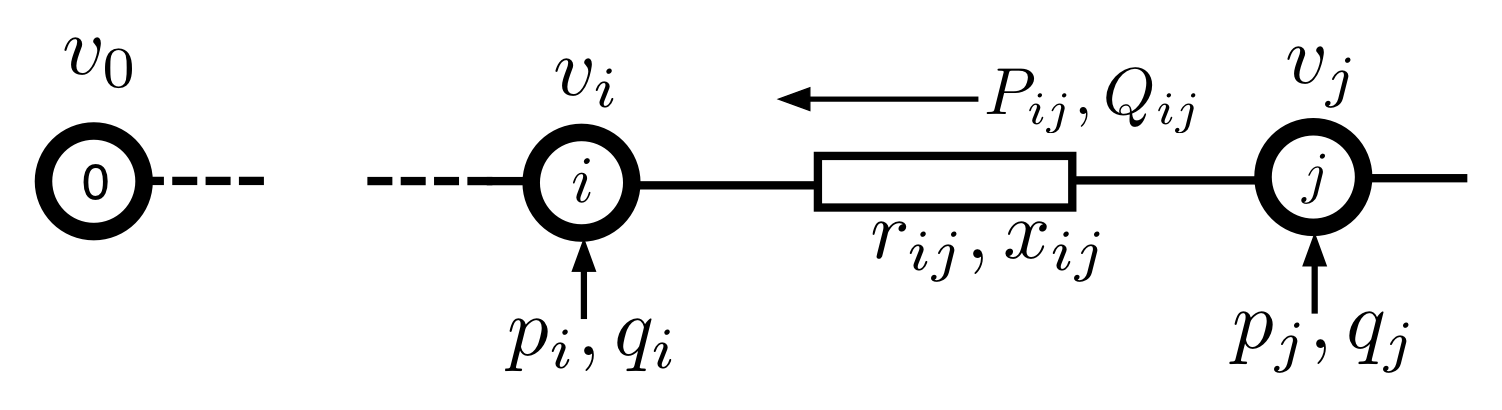}
\caption{\label{fig:radial_network} Nomenclature for a radial distribution network~\cite{heidari2017non}.}
\end{figure}

Consider a balanced, radial distribution network, shown in Fig.~\ref{fig:radial_network}, as an undirected graph $\mathcal{G}=\{\mathcal{N}\cup\{0\}, \mathcal{L}\}$ consisting of a set of $N+1$ nodes with $\mathcal{N}:=\{1, \hdots, N\}$ and a set of $N$ branches $\mathcal{L}:=\{1, \hdots, N\} \subseteq \mathcal{N}\times \mathcal{N}$, such that $(i,j)\in \mathcal{L}$, if nodes $i,j$ are connected. 
Node $0$ is assumed to be the substation node with a fixed voltage $V_0$. Let $B\in \mathbb{R}^{(N+1)\times N}$ be the \textit{incidence matrix} of $\mathcal{G}$ relating the branches in $\mathcal{L}$ to the nodes in $\mathcal{N}\cup \{0\}$, such that the $(i,k)$-th entry of $B$ is $1$, if the $i$-th node is connected to the $k$-th branch and, otherwise, $0$. Without loss of generality, $B$ can be organized to form an upper-triangular matrix. If $V_i$ is the voltage phasors at node $i$ and $I_{ij}$ is the current phasor in branch $(i,j)\in \mathcal{L}$, then define $v_i:=|V_i|^2$ and $l_{ij}:=|I_{ij}|^2$. Let $P_{ij}$ ($Q_{ij}$) be the active (reactive) power flow from node $j$ to $i$, denote $p_j$ ($q_j$) the active (reactive) power injections into node $j$, and let $r_{ij}$ ($x_{ij}$) be the resistance (reactance) of branch $(i,j)\in \mathcal{L}$, which means that the branch impedance is given by $z_{ij}:=r_{ij}+\mathbf{i}x_{ij}$. Then, for a radial network, the relation between nodal voltages and power flows is given by the \textit{DistFlow} equations $\forall (i,j)\in \mathcal{L}$:
\begin{subequations}\label{eq:dist_flow}
\begin{align}
v_j=&v_i+2r_{ij}P_{ij}+2x_{ij}Q_{ij}-|z_{ij}|^2l_{ij} \label{eq:volt_rel}\\
P_{ij}=&p_j+\sum_{h:h\rightarrow j}(P_{jh}-r_{jh}l_{jh}) \label{eq:real_power_rel}\\
Q_{ij}=&q_j+\sum_{h:h\rightarrow j}(Q_{jh}-x_{jh}l_{jh}) \label{eq:reac_power_rel}\\
l_{ij}(P_{ij},Q_{ij},v_j)  =& \frac{P_{ij}^2+Q_{ij}^2}{v_j}, \label{eq:curr_rel}
\end{align}
\end{subequations}
\textcolor{black}{where nodal power injections are $p_j:=p_{\text{g},j}-P_{\text{L},j}$ and $q_j:=q_{\text{g},j}-Q_{\text{L},j}$ with $p_{\text{g},j}$ ($q_{\text{g},j}$) as the controllable active (reactive) injections and $P_{\text{L},j}$ ($Q_{\text{L},j}$) is the uncontrollable active (reactive) demand.  The controllable injections include solar PV and flexible demand.}

The goal of this work is to maximize the range of active power injections, $p_{\text{g}}$, from a given feasible operating point with $p_{\text{g},j}=0,\, q_{\text{g},j}=0 \ \forall j\in \mathcal{N}$, such that all  voltages $v_j$ and currents $l_{ij}$ are within their respective limits (i.e., $v_j \in [\underline{v}_j,\overline{v}_j] \,\, \forall j\in \mathcal{N}$ and $l_{ij} \in [\underline{l}_{ij},\overline{l}_{ij}] \,\, \forall (i,j)\in \mathcal{L}$).
However, finding such a range is challenging due to the non-linear nature of~\eqref{eq:curr_rel}. For clarity, we provide definitions of the following key terms used in the manuscript.
\begin{definition}[AC Admissibility]
A solution, $p_{\text{g}}^\ast$, of a convex OPF problem is AC admissible, if the solution applied to~\eqref{eq:dist_flow} yields corresponding voltages and branch currents within their respective limits. 
\end{definition}
\begin{definition}[Nodal capacity]
Nodal capacity is the range of AC admissible active power injections $\Delta p_{\text{g},j} := [p_{\text{g},j}^-, p_{\text{g},j}^+]\subset \mathbb{R} \enspace \forall j\in \mathcal{N}$ with $p_{\text{g},j}^-\le 0$ and $p_{\text{g},j}^+\ge 0$. {\color{black} That is, the hyper-rectangle defined by the Cartesian product of intervals $\Delta p_{\text{g},j}$ contains all AC admissible (net) injections, i.e.,  all $p_{\text{g}} \in \Delta p_\text{g} := \times \Delta p_{\text{g},j} \subset \mathbb{R}^N$ are AC admissible, and characterizes the nodal capacities.}
\end{definition}

In the next section we use a simple 3-node system to motivate the need for analyzing nodal capacity in distribution systems. 

\subsection{Motivating example on nodal capacity}\label{sec:3-node_set}
\begin{figure}[t]
\centering
\includegraphics[width=0.3\textwidth]{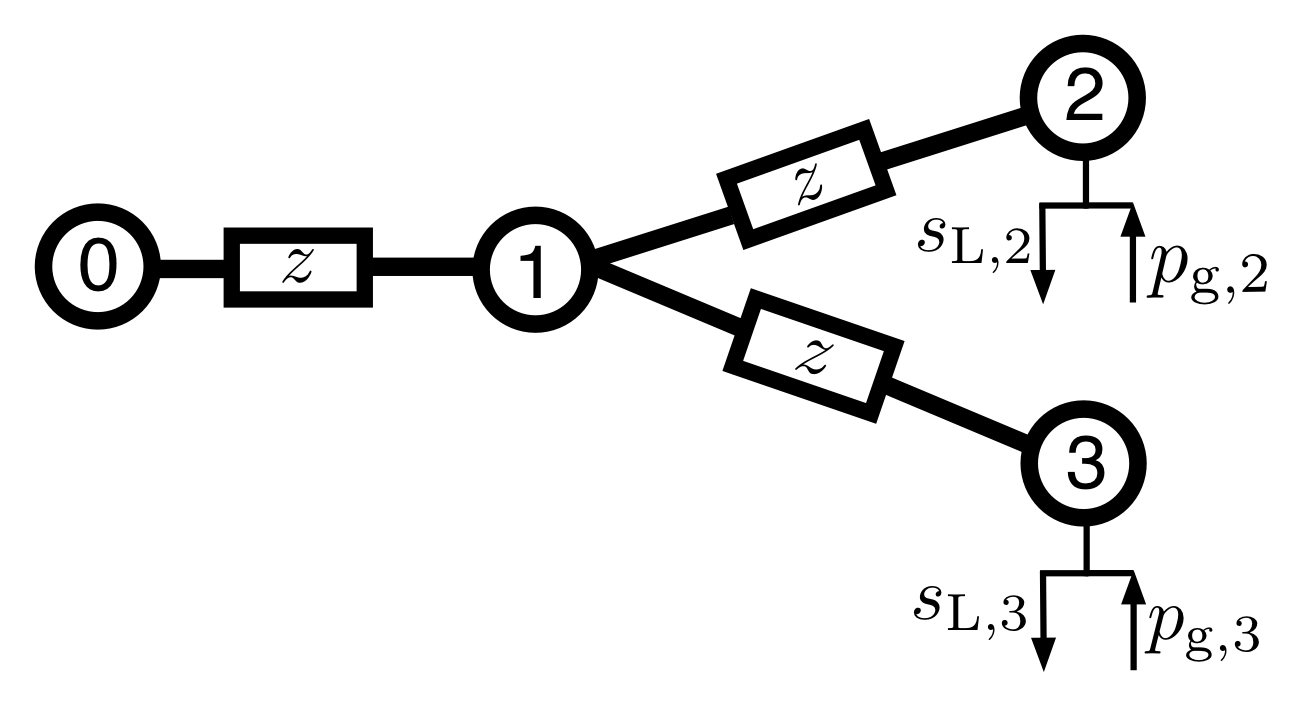}
\caption{\label{fig:3_node_model} The 3-node network used as a motivating example.}
\end{figure}

\begin{figure}[t]
\includegraphics[width=0.45\textwidth]{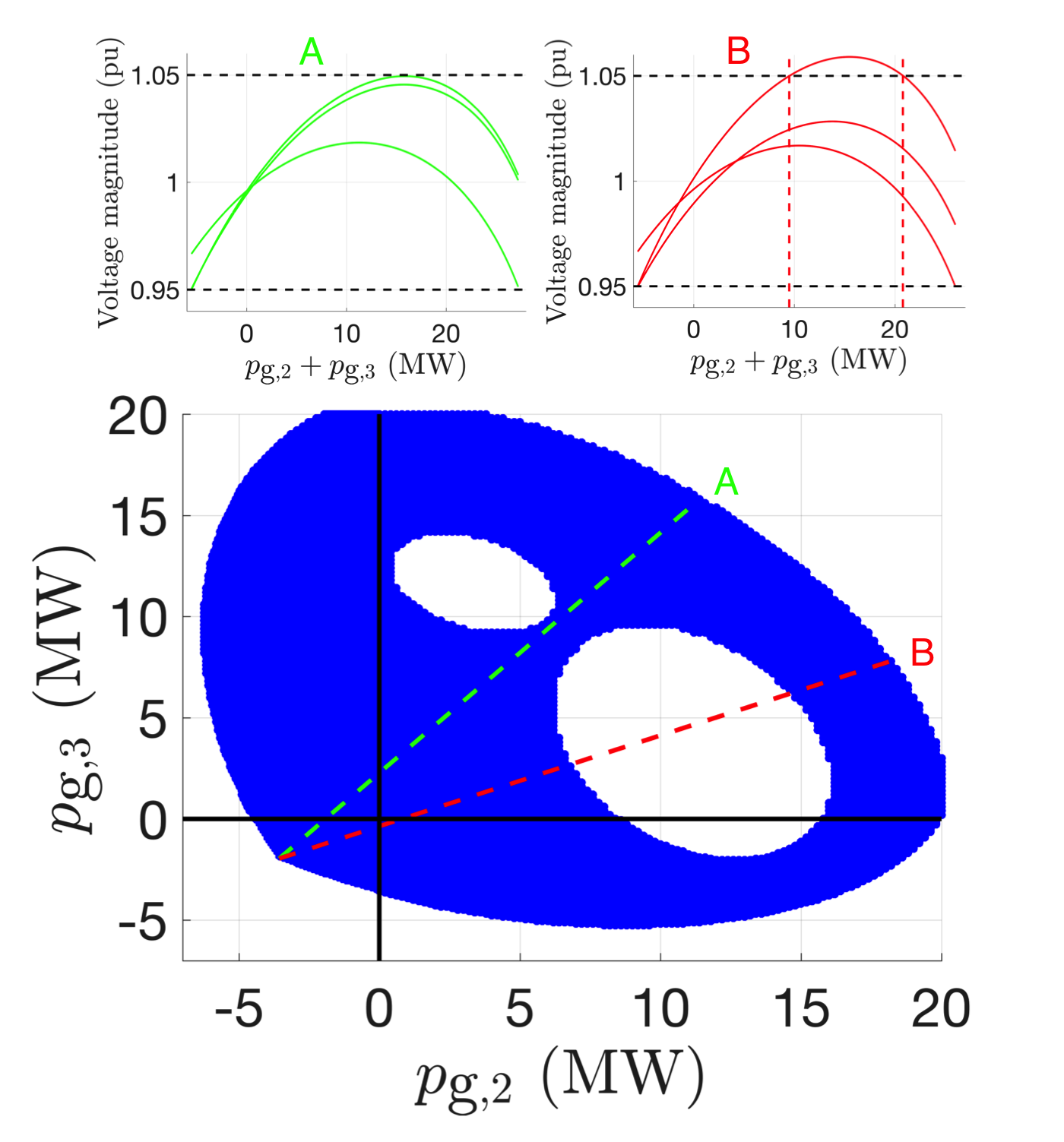}
\caption{\label{fig:3_node_sweep} AC admissibility for 3-node motivating example. (Top) Voltage profiles from sweeping $(p_{\text{g},2}, p_{\text{g},3})$ along admissible trajectory~A and inadmissible trajectory~B. (Bottom) The set of admissible injections is non-convex with trajectories A and~B showing admissible (green) and inadmissible (red) dispatch, respectively.}
\end{figure}

Consider the 3-node system shown in Fig.~\ref{fig:3_node_model}. Each branch of the system has an impedance of $z=0.55+\mathbf{i}1.33$pu. Node~$2$ has a load injection $s_{\text{L},2}=-0.02+\mathbf{i}0.005$pu and node~$3$ has a load injection $s_{\text{L},3}=-0.015+\mathbf{i}0.001$pu.
Flexible resources $p_{\text{g},2}$ and $p_{\text{g},3}$ are located at nodes $2$ and $3$. Only the active power at nodes $2$ and $3$ is controllable\footnote{For this manuscript, any mechanical devices such as tap-changers, capacitor banks and switches are assumed to be fixed at their nominal values and are not part of the optimization problem.}. 
By varying $p_{\text{g},2}$ and $p_{\text{g},3}$ and obtaining the corresponding AC power flows via Matpower~\cite{zimmerman2011matpower}, we can mapt the set of AC admissible injections as shown in Fig.~\ref{fig:3_node_sweep}(Bottom) with voltage limits (i.e., [0.95,1.05]). As can be seen from Fig.~\ref{fig:3_node_sweep}, the admissible set is non-convex and contains ``holes''. Hence, it is important when dispatching $p_{\text{g},2}$ and $p_{\text{g},3}$, to choose the right trajectory in order to maintain AC admissibility.
Fig.~\ref{fig:3_node_sweep} shows two different dispatch trajectories and the corresponding voltage profiles. Trajectory~A (green) is contained within the admissible set and, as a result, the voltage profile shown in green in the top of Fig.~\ref{fig:3_node_sweep}, satisfies the voltage limits. On the other hand, Trajectory~B (red) is not contained in the admissible set, which manifests itself as a voltage limit violation in the red curve in the top of Fig.~\ref{fig:3_node_sweep}. Even though trajectory~A is AC admissible it requires $p_{\text{g},2}$ and $p_{\text{g},3}$ to be coordinated (i.e., stay on the trajectory) to ensure admissibility, so they cannot be manipulated independently. This means that any changes in either requires a change in the other and, thus, end-point A in Fig.~\ref{fig:3_node_sweep}(Top) is not a nodal capacity. \textcolor{black}{ This is because a nodal capacity requires that nodal injections be manipulated independently.} \textcolor{black}{Note that due to the network coupling in the AC power flow equations, the nodal capacities are inherently coupled and cannot be solved independently of each other.}
This simple example shows the need to develop tools that determine nodal capacities for any radial, balanced network. Towards that objective, the next section develops a convex inner approximation of the nonlinear \textit{DistFlow} formulation in~\eqref{eq:dist_flow}. 

\subsection{Convex Inner Approximation Preliminaries}\label{sec:bounds}
In this section, we first present a compact matrix representation of the linear components~\eqref{eq:volt_rel}-\eqref{eq:reac_power_rel}.
Then, we bound the nonlinear branch current terms, $l_{ij}(P_{ij},Q_{ij},v_j)$ in~\eqref{eq:curr_rel}, by a convex envelope, which begets a convex inner approximation of~\eqref{eq:dist_flow}.


First, define vectors $P:=[P_{ij}]_{(i,j)\in \mathcal{L}}\in \mathbb{R}^N$, $Q:=[Q_{ij}]_{(i,j)\in \mathcal{L}}\in \mathbb{R}^N$, $V:=[v_i]_{i\in \mathcal{N}}\in \mathbb{R}^N$, $p:=[p_i]_{i\in \mathcal{N}}\in \mathbb{R}^N$, $p_{\text{g}}:=[p_{\text{g},i}]_{i\in \mathcal{N}}\in \mathbb{R}^N$, $P_{\text{L}}:=[P_{\text{L},i}]_{i \in \mathcal{N}} \in \mathbb{R}^N$, $q:=[q_i]_{i \in \mathcal{N}}\in \mathbb{R}^N$, $Q_{\text{L}}:=[Q_{\text{L},i}]_{i\in \mathcal{N}} \in \mathbb{R}^N$, and  $l:=[l_{ij}]_{(i,j)\in \mathcal{L}}\in \mathbb{R}^N$ and matrices $R:=\text{diag}\{r_{ij}\}_{(i,j)\in \mathcal{L}}\in \mathbb{R}^{N\times N}$, $X:=\text{diag}\{x_{ij}\}_{(i,j)\in \mathcal{L}}\in \mathbb{R}^{N\times N}$, $Z^2:=\text{diag}\{z_{ij}^2\}_{(i,j)\in \mathcal{L}}\in \mathbb{R}^{N\times N}$, and $A:=[0_N \quad I_N]B-I_N$, where $I_N$ is the $N\times N$ identity matrix and $0_N \in \mathbb{R}^N$. Then, directly applying~\cite{heidari2017non}, we get the following expression for $P$, $Q$ and $V$:
\begin{align}
    V=&v_{\text{0}}\mathbf{1}_N+M_{\text{p}}p+M_{\text{q}}q-Hl,\label{eq:final_volt_rel}\\
    P=&Cp-D_{\text{R}}l, \qquad Q=Cq-D_{\text{X}}l,\label{eq:P_relation}
\end{align}
where matrices $M_{\text{p}}:=2C^TRC$, $M_{\text{q}}:=2C^TXC$, $C:=(I_N-A)^{-1}$, $H:=C^T(2(RD_{\text{R}}+XD_{\text{X}})+Z^2)$, $D_{\text{R}}:=(I_N-A)^{-1}AR$, and $D_{\text{X}}:=(I_N-A)^{-1}AX$ describe the network topology and impedance parameters. Note that in the authors' previous work in ~\cite{nazir2019convex}, it is shown that the matrix $(I_N-A)$ is non-singular for radial and balanced distribution networks. Furthermore, the convex inner approximation in~\cite{nazir2019voltage} is valid only for purely inductive\footnote{The term inductive (capacitive) branch refers herein to a network branch whose reactance is inductive (capacitive), which means $x_{ij}\ge 0$ ($x_{ij}<0$).}, radial, and balanced networks.
In the current manuscript, we extend the convex formulation to any radial and balanced network, including those with mixed inductive and capacitive branches.

Clearly, \eqref{eq:final_volt_rel} and~\eqref{eq:P_relation} represent linear relationships between the nodal power injections, $(p,q)$, the branch power flows, $(P,Q)$, and node voltages $V$. 
However, setting $l=0$ and neglecting~\eqref{eq:curr_rel}, as done with the commonly used \textit{LinDist} approximation, can result in overestimating the nodal capacities~\cite{nazir2019convex}. Next, we present methods for bounding the nonlinearity $l_{ij}(P_{ij},Q_{ij},v_j)$ from above and below.

Based on the description of voltages in~\eqref{eq:final_volt_rel} and branch flows in~\eqref{eq:P_relation}, denote $l_{\text{lb}}$ and $l_{\text{ub}}$ as lower and upper bounds on $l$. Then, we can define proxy variables for the corresponding upper $(.)^+$ and lower $(.)^-$ bounds of $P$, $Q$ and $V$ as follows:
\begin{subequations}\label{eq:CIA_bounds}
\begin{align}
    P^+:=&Cp-D_{\text{R}}l_{\text{lb}} \label{eq:P_relation_1}\\
    P^-:=&Cp-D_{\text{R}}l_{\text{ub}} \label{eq:P_relation_2}\\
    Q^+:=&Cq-D_{\text{X}_+}l_{\text{lb}}-D_{\text{X}_-}l_{\text{ub}}\label{eq:Q_relation_1}\\
    Q^-:=&Cq-D_{\text{X}_+}l_{\text{ub}}-D_{\text{X}_-}l_{\text{lb}}\label{eq:Q_relation_2}\\
    V^+:=&v_{\text{0}}\mathbf{1}_n+M_{\text{p}}p+M_{\text{q}}q-H_+l_{\text{lb}}-H_-l_{\text{ub}}\label{eq:V_relation_1}\\
    V^-:=&v_{\text{0}}\mathbf{1}_n+M_{\text{p}}p+M_{\text{q}}q-H_+l_{\text{ub}}-H_{-} l_{\text{lb}}, \label{eq:V_relation_2}
\end{align}
\end{subequations}
where $D_{\text{X}_+}$ and $H_+$ include the non-negative elements of $D_{\text{X}}$ and $H$, respectively, and $D_{\text{X}_-}$ and $H_-$ are the corresponding negative elements. For example, if the network is purely inductive, then $D_{\text{X}_-}=H_-=0$ and the formulation reduces to the one presented in~\cite{nazir2019voltage}. These upper and lower bounds in~\eqref{eq:CIA_bounds} satisfy $P^-\le P \le P^+$, $Q^-\le Q\le Q^+$ and $V^-\le V\le V^+$.
Note that bounds $l_\text{lb}$ and $l_\text{ub}$ in~\eqref{eq:CIA_bounds} effectively allow us to neglect the nonlinear~\eqref{eq:curr_rel}. Thus, if we can find convex representations of these bounds, the corresponding OPF formulation will be a convex inner approximation. This is described next. 

     Equation~\eqref{eq:CIA_bounds} provides a linear formulation for bounding the AC power flow equations in terms of bounds $l_\text{lb}, l_\text{ub}$ and controllable injections. 
     This was first presented in~\cite{nazir2019voltage}, where bounds $l_\text{lb}, l_\text{ub}$ were derived based on a nominal operating point and used to maximize voltage margins with mechanical grid assets (e.g., LTCs and capacitor-banks). Next, we summarize the derivation of these bounds and leverage them to formulate a novel convex inner approximation of the AC OPF to determine the nodal capacities.
  
Based on nominal operating point $x^0_{ij}:=\text{col}\{P_{ij}^0, Q_{ij}^0, v_j^0\} \in\mathbb{R}^3$, the second-order approximation for~\eqref{eq:curr_rel} can be expressed as
\begin{align}\label{eq:T_exp}
    l_{ij} & \approx l_{ij}^0 + \mathbf{J}_{ij}^\top \mathbf{\delta}_{ij} +\frac{1}{2}\mathbf{\delta}_{ij}^\top \mathbf{H}_{\text{e},ij} \mathbf{\delta}_{ij},
\end{align}
{\color{black}where $l_{ij}^0 := l_{ij}(x_{ij}^0)$ are squared branch currents and $\mathbf{\delta}_{ij}$, Jacobian $\mathbf{J}_{ij}$, and Hessian $\mathbf{H}_{\text{e},ij}$ are defined below:
\begin{align}
\mathbf{\delta}_{ij}:=
    \begin{bmatrix} 
        P_{ij}-P_{ij}^0\\
        Q_{ij}-Q_{ij}^0\\ 
        v_j-v_j^0
    \end{bmatrix},
    ~
\mathbf{J}_{ij}:=\left. \begin{bmatrix}
\frac{\partial l_{ij}}{\partial P_{ij}}\\
\frac{\partial l_{ij}}{\partial Q_{ij}}\\
\frac{\partial l_{ij}}{\partial v_{j}}\\
\end{bmatrix} 
\right|_{x_{ij}^0}
=
\begin{bmatrix}
\frac{2P^0_{ij}}{v^0_j}\\
\frac{2Q^0_{ij}}{v^0_j}\\ 
-\frac{(P^0_{ij})^2+(Q^0_{ij})^2}{(v^0_j)^2}
\end{bmatrix}\\
\mathbf{H}_{\text{e},ij}:=
\begin{bmatrix} 
    \frac{2}{v^0_j} && 0 && \frac{-2P^0_{ij}}{(v^0_j)^2}\\
    0 && \frac{2}{v^0_j} && \frac{-2Q^0_{ij}}{(v^0_j)^2}\\
    \frac{-2P^0_{ij}}{(v^0_j)^2} && \frac{-2Q^0_{ij}}{(v^0_j)^2} && 2\frac{(P^0_{ij})^2+(Q^0_{ij})^2}{(v^0_j)^3}
    \end{bmatrix}.
\end{align}
}

\begin{figure}[h]
\centering
\includegraphics[width=0.4\textwidth]{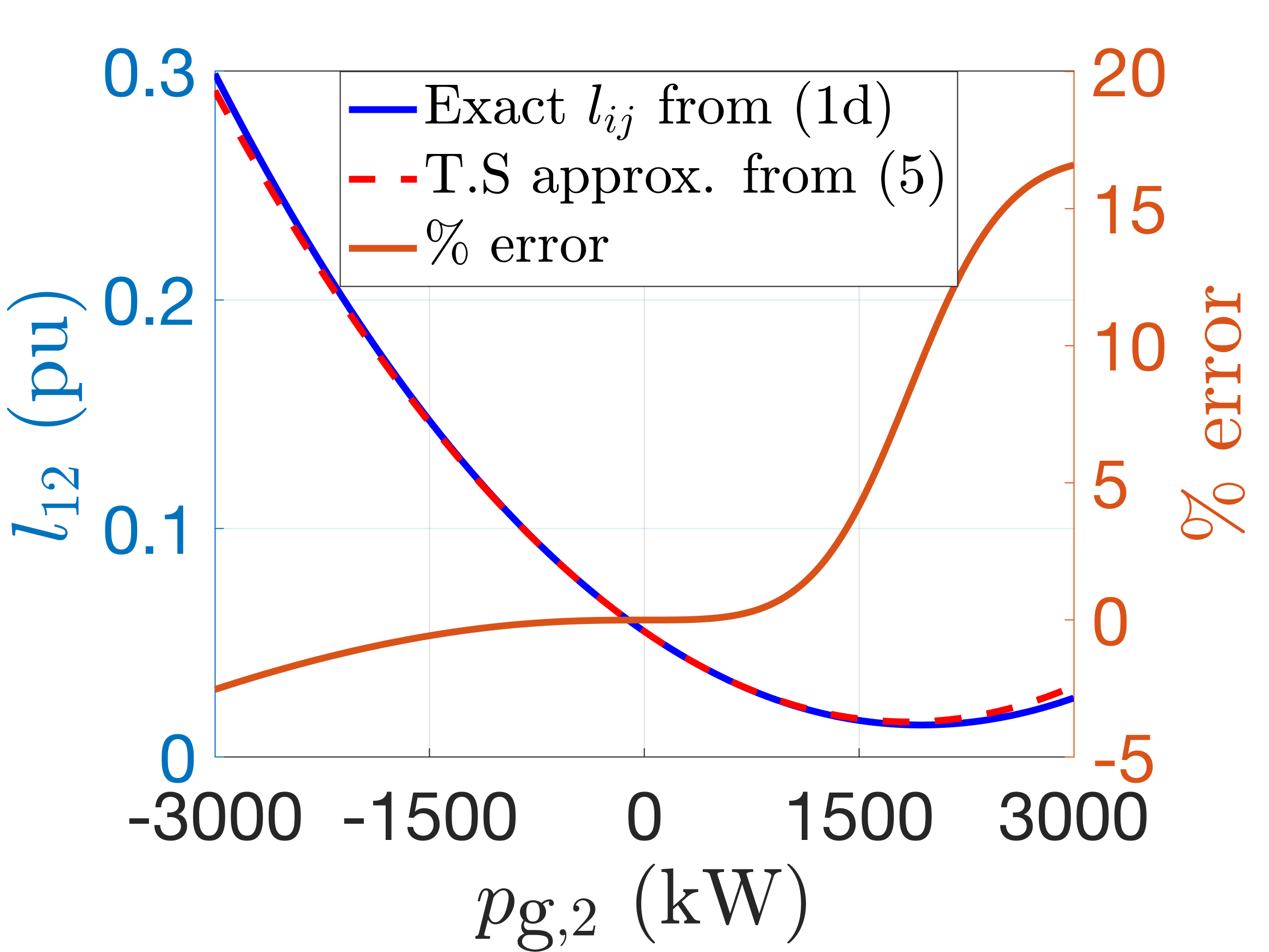}
\caption{\label{fig:TS_approx} {\color{black} Comparison of the second-order approximation from~\eqref{eq:T_exp} with the original (nonlinear) expression of $l$ in~\eqref{eq:curr_rel} for the IEEE 13-node feeder. The maximum absolute approximation error is found at branch $l_{12}$ when the active net-injection at node~$2$ ($p_{\text{g},2}$) is $-3000$kW and is $0.0074$pu (and $\le 2.5\%$) and relative error is less than 5\% for $p_{\text{g},2} \in [-3000,1600]$ kW.}}
\end{figure}

Fig.~\ref{fig:TS_approx} compares the accuracy of the second-order approximation in~\eqref{eq:T_exp} with the non-linear expression from~\eqref{eq:curr_rel}. Note that, for the IEEE-13 node network, the worst-case approximation error for $l_{ij}$ is less than $0.008$pu over the wide range of net injections $[-3000,3000]$kW. Thus, we can assume that the expression in~\eqref{eq:T_exp} is sufficiently accurate and omit the higher-order terms.
 
 Furthermore,~\cite{nazir2019voltage} shows that $\mathbf{H}_{\text{e},ij}$ is positive semi-definite, which, together with~\eqref{eq:T_exp}, means that lower and upper bounds of $l_{ij}$ for all $(i,j)\in \mathcal{L}$ are given by:
 \begin{align}
    l_{ij}=|l_{ij}| & \approx |l_{ij}^0 + \mathbf{J}_{ij}^\top\mathbf{\delta}_{ij}+\frac{1}{2}\mathbf{\delta}_{ij}^\top \mathbf{H}_{\text{e},ij} \mathbf{\delta}_{ij}| \\
    & \le |l_{ij}^0| + |\mathbf{J}_{ij}^\top\mathbf{\delta}_{ij}|+|\frac{1}{2}\mathbf{\delta}_{ij}^\top \mathbf{H}_{\text{e},ij} \mathbf{\delta}_{ij}| \\
    & \le l_{ij}^0 + \max\{2|\mathbf{J}_{ij}^\top\mathbf{\delta}_{ij}|,|\mathbf{\delta}_{ij}^\top \mathbf{H}_{\text{e},ij} \mathbf{\delta}_{ij}|\}\label{eq:l_upper_1}\\
    \implies  l_{ij} &
     \le l_{ij}^0 + \max\{2|\mathbf{J}_{ij+}^\top\mathbf{\delta}_{ij}^++\mathbf{J}_{ij-}^\top\mathbf{\delta}_{ij}^-|,\mathbf{\psi}_{ij}\} \le  l_{\text{ub},ij}, \label{eq:l_upper}\\
     \text{and } l_{ij} & \ge l_{ij}^0 + \mathbf{J}_{ij+}^\top\mathbf{\delta}_{ij}^-+\mathbf{J}_{ij-}^\top\mathbf{\delta}_{ij}^+ =: l_{\text{lb},ij}, \label{eq:l_lower}
\end{align}
where $\mathbf{J}_{ij+}$ and $\mathbf{J}_{ij-}$ includes the positive and negative elements of $\mathbf{J}_{ij}$, $\mathbf{\delta}_{ij}^+ := \mathbf{\delta}_{ij}(P_{ij}^+, Q_{ij}^+, v_j^+,x_{ij}^0)$ and $\mathbf{\delta}_{ij}^-:=\mathbf{\delta}_{ij}(P_{ij}^-, Q_{ij}^-, v_j^-,x_{ij}^0)$, and $\mathbf{\psi}_{ij}:=\max\{(\mathbf{\delta}_{ij}^{+,-})^\top \mathbf{H}_{\text{e},ij} (\mathbf{\delta}_{ij}^{+,-})\}$, which represents the largest of eight possible combinations of $P/Q/v$ terms in $\mathbf{\delta}_{ij}$ with mixed $+,-$ superscripts.
Note that from~\eqref{eq:l_lower}, the lower bound $l_{\text{lb},ij}$ may become negative, however, we know from the physics that $l_{ij}\ge 0$, which means the $l_{\text{lb},ij}$ may be conservative.
{\color{black}Also, from the expression of $\mathbf{J}_{ij}$, it can be seen that if $P_{ij}^0,\enspace Q_{ij}^0\approx 0$, then $l_{\text{lb},ij}=l_{ij}^0\approx0$ and so the lower bound is fixed for this operating point. This results in a conservative estimate of the lower bound, which can result in undersizing nodal hosting capacity.}
To alleviate this shortcoming, Algorithm~\ref{alg1} in Section~\ref{sec:Iter_alg} presents an iterative approach that improves the nodal capacity.
Thus, with~\eqref{eq:CIA_bounds},~\eqref{eq:l_upper}, and~\eqref{eq:l_lower}, we have a convex inner approximation of~\eqref{eq:dist_flow} that can be used to determine the nodal capacities. \textcolor{black}{Furthermore, since the Jacobian and Hessian are calculated separately for each branch in the network, the size of the matrices does not grow with the size of the system. This enables the extension of this approach to large scale systems.} 


\subsection{Optimizing DER nodal capacity}\label{sec:CIA_form} 
The bounds from~\eqref{eq:l_upper} and~\eqref{eq:l_lower} allow us to omit~\eqref{eq:curr_rel} entirely and replace the original variables $P$, $Q$, and $V$ with corresponding proxy variables that serve as upper and lower bounds $(.)^+$ and $(.)^-$ in~\eqref{eq:CIA_bounds}. Since $(.)^+$ and $(.)^-$ are outer approximations, using them in an OPF formulation results in a feasible set that is contained in the original, non-convex AC OPF. This means that the following represents a convex inner approximation and can be used to determine nodal capacities:
\begin{subequations}\label{eq:P1}
\begin{align}
\text{(P1)}  \quad p_{\text{g}}^+(p_{\text{g}}^-)=\arg\min_{p_{\text{g},i},q_{\text{g},i}} \  \sum_{i=1}^{N}f_i(p_{\text{g},i})&\\
\text{s.t.}  \quad  \eqref{eq:P_relation_1}-\eqref{eq:V_relation_2},\eqref{eq:l_upper},& \eqref{eq:l_lower} \\
    p=p_{\text{g}}-P_{\text{L}} \quad q=  q_{\text{g}}-Q_{\text{L}}, & \\
    \underline{V}\le     V^-(p,q) \quad V^+(p,q)\le &\overline{V}   \label{eq:P1_V}\\
    l_{\text{ub}}\le  \overline{l} \quad \underline{q_{\text{g}}}\le q_{\text{g}}\le \overline{q_{\text{g}}}. \label{eq:P1_lmax}
\end{align}
\end{subequations}
 In (P1), ~\eqref{eq:P1_V} ensures that any feasible dispatch $p_\text{g}$ satisfies the nodal voltages in the original AC OPF based on~\eqref{eq:dist_flow}. {\color{black} Also, the current limit constraint in~\eqref{eq:P1_lmax} ensures that the upper limit of the branch current ($l_{\text{ub}}$) is within the given thermal line current limits.} To determine the hyper-rectangle $\Delta p_\text{g}$,
 we must solve (P1) once for the lower ranges, $p_{\text{g}}^-$, and once for the upper ranges, $p_{\text{g}}^+$. Thus, the objective function components, $f_i(p_{\text{g},i})$, must be designed to engender $p_{\text{g},i}^-$ and $p_{\text{g},i}^+$. For example, to compute $p_{\text{g},i}^-$, we can choose $f_i(p_{\text{g},i}):=\alpha_i p_{\text{g},i}\,$ and, for $p_{\text{g},i}^+$, we can designate $f_i(p_{\text{g},i}):=-\alpha_i p_{\text{g},i}$, where $\alpha_i\ge 0$ is the relative priority  of nodal capacity at node $i$. Clearly, the choice of objective function determines how flexibility is allocated over the network, e.g., choosing objective function such as $\pm \alpha_i\log(p_{\text{g},i})$ can result in a different allocation of nodal capacity over the network as compared with $\pm \alpha_i p_{\text{g},i}$. The design of the objective function represents an interesting future extension into energy policy and incentive mechanism and rate design~\cite{perez2017regulatory}.

While (P1) ensures AC admissibility at the nodal capacity values, it is natural to consider what happens when the nodal flexibility is below the rated capacity. That is, are all injections within the hosting capacity range guaranteed to be admissible across all the nodes? The next section answers this question by providing analytical guarantees of admissibility for the nodal hosting capacity, within the hyper-rectangle $\Delta p_{\text{g}}$, and then presents an iterative algorithm to successively improve $\Delta p_{\text{g}}$.

\section{Analysis of convex inner approximation}\label{sec:Feas_guar}
Next, we analyze (P1) and prove that all (net) injections $p_{\text{g}} \in \Delta p_{\text{g}}$ are AC admissible. 

\subsection{Admissibility guarantees}
Due to page limitations, AC admissibility guarantees consider only nodal voltage limits, however, the case of admissibility under branch flow limits is similar.
\begin{theorem}\label{theorem1}
Under conditions C1)~$\frac{\partial V^+}{\partial p_{\text{g},i}}\ge 0$, C2)~$\frac{\partial V^-}{\partial p_{\text{g},i}}\ge 0$, $\forall i~\in~\mathcal{N}$, if $\Delta p_{\text{g}}$ is obtained via nodal capacities from (P1), then $\forall p_{\text{g}}\in \Delta p_{\text{g}}$ and $p(p_{\text{g}}) = p_{\text{g}} - P_L$, we have 
$$\underline{V} \le V^-(p) \le V(p)\le V^+(p)\le \overline{V},$$
 where $V(p)$ represents the actual nodal voltages from~\eqref{eq:dist_flow} resulting from (net) injections $p$.
\end{theorem}
\begin{proof}
Consider two cases: $0 \le p_{\text{g}} \le p_{\text{g}}^+$ (Case~1) and $p_{\text{g}}^- \le p_{\text{g}} \le 0$ (Case~2).

\noindent
\textit{Proof of Case~1:} Using~\eqref{eq:V_relation_1} at $p_{\text{g}}^+$ yields:
\begin{align}
    V^+(p^+)=v_0\mathbf{1}_n+M_{\text{p}}p^+ + M_{\text{q}}q^+-H_+l_{\text{lb}}-H_-l_{\text{ub}}\le \overline{V}
\end{align} 
where $p^+=p_{\text{g}}^+-P_{\text{L}}$ and $q^+=q_{\text{g}}^+-Q_{\text{L}}$. Now, consider any $p_{\text{g}}\in \Delta p_{\text{g}}$ such that $0 \le p_{\text{g}}\le p_{\text{g}}^+$ and using C1, then
\begin{align}\label{eq:V_nlmodel}
V^+(p)=v_0\textbf{1}_n+M_{\text{p}}p + M_{\text{q}}q^+-H_+l_{\text{lb}}(p)-H_-l_{\text{ub}}(p)\le\overline{V}
\end{align}
where $p=p_{\text{g}}-P_{\text{L}}$. The actual voltage according to~\eqref{eq:final_volt_rel} at $p$ is
\begin{align}\label{eq:V_linmodel}
    V(p)=v_0\textbf{1}_n+M_{\text{p}}p+M_{\text{q}}q^+-H_+l(p)-H_-l(p)
\end{align}
Then, subtracting \eqref{eq:V_nlmodel} from \eqref{eq:V_linmodel} gives:
\begin{align}
    V^+(p)-V(p)=H_+\left(l(p)-l_{\text{lb}}(p)\right)+H_-\left(l(p)-l_{\text{ub}}(p)\right)
\end{align}
Using~\eqref{eq:l_lower} and~\eqref{eq:l_upper} we get, $l_{\text{lb}}(p)\le l(p)\le l_{\text{ub}}(p)$ and that $V^+(p)-V(p)\ge 0 \implies V(p)\le V^+(p)\le \overline{V}$.

\noindent
\textit{Proof of Case~2:} 
Using~\eqref{eq:V_relation_2} at $p_{\text{g}}^-$ yields:
\begin{align}
    V^-(p^-)=v_0\mathbf{1}_n+M_{\text{p}}p^-+M_{\text{q}}q^--H_+l_{\text{ub}}-H_-l_{\text{lb}}\ge \underline{V}
\end{align} 
where $p^-=p_{\text{g}}^--P_{\text{L}}$ and $q^-=q_{\text{g}}^--Q_{\text{L}}$. Now, consider any $p_{\text{g}}\in \Delta p_{\text{g}}$ such that $p_{\text{g}}^- \le p_{\text{g}} \le 0$ and C2, then
\begin{align}\label{eq:V_negmodel}
V^-(p)=v_0\textbf{1}_n+M_{\text{p}}p+M_{\text{q}}q^--H_+l_{\text{ub}}(p)-H_-l_{\text{lb}}(p)\ge\underline{V}
\end{align}
where $p=p_{\text{g}}-P_{\text{L}}$. The actual voltage according to~\eqref{eq:final_volt_rel} at $p$ is
\begin{align}\label{eq:V_linmodel2}
    V(p)=v_0\textbf{1}_n+M_{\text{p}}p+M_{\text{q}}q^--H_+l(p)-H_-l(p)
\end{align}
Then, subtracting \eqref{eq:V_negmodel} from \eqref{eq:V_linmodel2} gives:
\begin{align}
     V^-(p)-V(p)=H_+\left(l(p)-l_{\text{ub}}(p)\right)+H_-\left(l(p)-l_{\text{lb}}(p)\right)
\end{align}
Using~\eqref{eq:l_lower} and~\eqref{eq:l_upper} we get, $l_{\text{lb}}(p)\le l(p)\le l_{\text{ub}}(p)$, and that $V^-(p)-V(p)\le 0 \implies V(p)\ge V^-(p)\ge \underline{V}$.
Combining Case~1 and Case~2 completes the proof.
\end{proof}
	\textcolor{black}{Note that conditions C1 and C2 are not trivial condition and may fail to hold under certain operating conditions. This can be seen from~\eqref{eq:V_relation_1} and~\eqref{eq:V_relation_2}, as $V^+$ and $V^-$ depend on the upper bound  $l_{\text{ub}}$ of the loss term $l$. This upper bound in turn is based on the expression in~\eqref{eq:l_upper} and can either increase or decrease with $p$. This is further illustrated in Fig.~\ref{fig:3_node_sweep}, which shows that the voltage magnitude can decrease with increasing power injections.}
		\textcolor{black}{In addition, conditions C1 and C2 only guarantee satisfaction of network constraints for injections in $\Delta p_{\text{g}}$ but do not guarantees existence of AC power flow solutions, which is provided in Section~III.D.}
Finally, Theorem~\ref{theorem1} significantly improves over the result provided in~\cite{nazir2019voltage}, since it guarantees that the entire set, $\Delta p_{\text{g}}$, is AC admissible rather than just the solutions, $p_{\text{g}}^+$ and $p_{\text{g}}^-$. Importantly, this is \textit{exactly} why $\Delta p_{\text{g}}$ can be used to characterize the nodal hosting capacity.
As with any convex inner approximation, the results can be conservative. Thus, in the next section, a new iterative algorithm is presented that successively increases the nodal capacity.

\subsection{Iterative algorithm for nodal capacity improvement}\label{sec:Iter_alg}
The bounds $l_\text{lb}, l_\text{ub}$ obtained in Section~\ref{sec:bounds} can be conservative depending upon the nominal operating point, $x^0$.
Thus, when we solve (P1) to determine $p_{\text{g}}^+$ and $p_{\text{g}}^-$, the nodal capacities can be significantly underestimated. This is especially the case when $P_{ij}^0, Q_{ij}^0\approx0$, which yields a Jacobian close to zero and the first-order estimate of $l_{\text{lb},ij}$ also close to $l_{ij}^0\approx0$ per~\eqref{eq:l_lower}. 
Algorithm~1 mitigates the conservativeness by successively augmenting the set $\Delta p_{\text{g}}$ via updating the operating point, Jacobian, and Hessian to reflect known nodal capacities. This approach is similar to the so-called convex-concave procedure~\cite{rosenlicht1986introduction}. Algorithm~\ref{alg1} outlines the steps involved in the proposed scheme.
\begin{algorithm}
\caption{\label{alg1} Successive enhancement of DER nodal capacity $\Delta p_{\text{g}}$ (unity power factor case)}
\SetAlgoLined
\KwResult{Admissible set of injections $\Delta p_{\text{g}} := \times [p^-_{\text{g},i},p^+_{\text{g},i}]$}
\textbf{Input:} $P_L, Q_L\in \mathbb{R}^N$, convex $f_i(p_{\text{g},i})\,\, \forall i\in \mathcal{N}$, and $\epsilon>0$\\
 Run Load flow w/ $P_L,Q_L,p_{\text{g}}(0)=0_N$ $\Rightarrow \mathbf{J}(0),\mathbf{H}_{\text{e}}(0)$\\
 \For{$m=1:2$}{
 \eIf{$m=1$}{
 $p_{\text{g},i} \rightarrow p_{\text{g},i}^+$, Cond($i$)$\rightarrow$ Check $\frac{\partial V^+}{\partial p_{\text{g},i}}\ge 0 \quad \forall i \in \mathcal{N}$}
     { 
        $p_{\text{g},i} \rightarrow p_{\text{g},i}^-$, Cond($i$)$\rightarrow$ Check $\frac{\partial V^-}{\partial p_{\text{g},i}}\ge 0 \quad \forall i \in \mathcal{N}$
        }
 Initialize $k=1$,  $error(0)=\infty$\\
  \While{$\exists i$, s.t. Cond($i$) holds $\land$  $error(k-1) > \epsilon$}{
  \For{$i=1:N$}{
  \If{ Cond(i) does not hold}{Set $p_{\text{g},i}(k)=p_{\text{g},i}(k-1)$}
  }
        Solve (P1) $\Rightarrow p_{\text{g},i}(k), f_i(p_{\text{g},i}(k)), \ \forall i \in \mathcal{N}$\\
        Run load flow w/ $P_L-p_{\text{g}}(k), Q_L \Rightarrow \mathbf{J}(k), \mathbf{H}_{\text{e}}(k)$\\
        Update Cond($i$) $\forall i \in \mathcal{N}$\\
        Update error: \\ 
        {\small $error(k)=\max_{i\in\mathcal{N}} \left|f_i(p_{\text{g},i}(k))-f_i(p_{\text{g},i}(k-1))\right|$}\\
        $k \rightarrow k+1$
        }
        }
\end{algorithm}

Next, we illustrate how the nodal capacity improves with Algorithm~\ref{alg1} for the motivating example in Fig.~\ref{fig:3_node_model}. 
Note that reactive power net injections, $q_{\text{g},i}$, are decision variables in (P1), however, in the proceeding analysis and simulations, we set $q_{\text{g},i}=0$. Later in Section~\ref{sec:Reac_power}, we analyze the role of reactive power strategies to further augment nodal capacities. For the sake of simplicity, we neglect the branch limit constraint~\eqref{eq:P1_lmax} in (P1) and assume an oversized substation transformer, which is a common practice in the US. The focus is on voltage because that is often the primary concern of utilities in the US~\cite{HOROWITZ2020106222}. However, the formulation in (P1) and the analysis therein hold for branch limit constraints as well. 

\begin{figure}[t]
\centering
\includegraphics[width=0.33\textwidth]{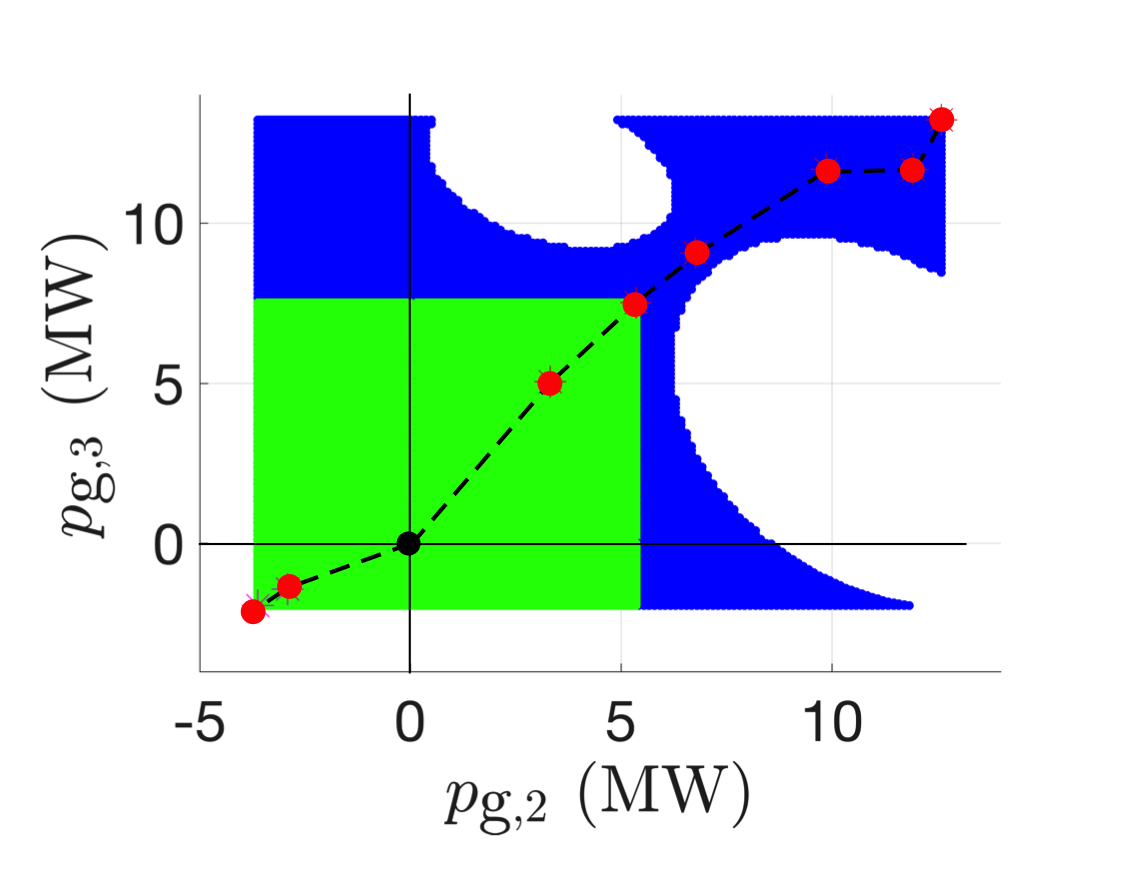}
\caption{\label{fig:3_node_sweep_mon} The set of admissible injections for the 3-node network is non-convex (blue). Algorithm~\ref{alg1} can find maximal admissible injections via iterations (red dots), but  monotonicity conditions $C1$ and $C2$ in Theorem~\ref{theorem1} are important to define the convex inner approximation (green), which yields the nodal capacity, $\Delta p_{\text{g}}$.}
\end{figure}

{\color{black} The importance of conditions C1 and C2 in Algorithm~\ref{alg1} in determining $\Delta p_{\text{g}}$ is illustrated in Fig.~\ref{fig:3_node_sweep_mon}. Specifically, without C1 and C2 (i.e., omitting lines 12-14 from Algorithm~\ref{alg1}), we get a piece-wise linear dispatch trajectory (red dots and dashed black lines), which represents a large nodal dispatch range, but does not permit independent dispatch from different nodes and is, thus, not a valid nodal capacity. 
However, when the conditions C1 and C2 are explicitly considered then Algorithm~\ref{alg1} yields the green region in Fig.~\ref{fig:3_node_sweep_mon}, which is a convex inner approximation of the blue set (of AC admissible dispatches).} 
From the non-convex nature of the blue set, it is clear that operating beyond the green set may require coordination between different nodes in order to ensure AC admissibility. 
Furthermore, since the green set represents $\Delta p_{\text{g}}$ and is a hyper-rectangle, no coordination between nodes is necessary in order to guarantee AC admissibility. 
In summary, Algorithm~\ref{alg1} can iteratively solve (P1) to augment $\Delta p_{\text{g}}$ to achieve a larger nodal capacity, which reduces conservativeness of the CIA-based approach. 

Next, we present Case Study~1, which employs Algorithm~\ref{alg1} to determine the solar PV hosting capacity for a distribution network.

\subsubsection*{\textbf{Case study 1}}
\begin{figure}[t]
\centering
\includegraphics[width=0.28\textwidth]{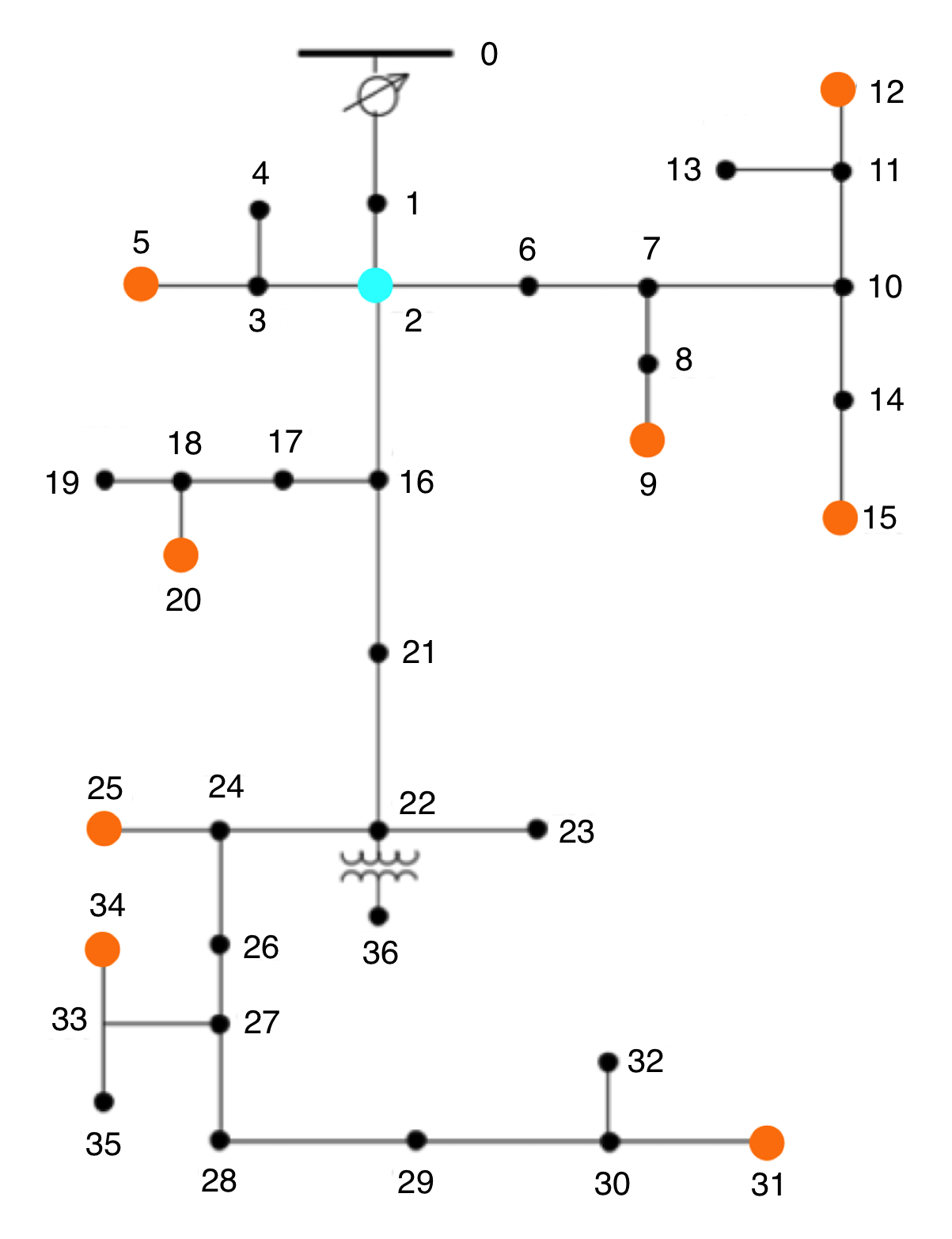}
\caption{\label{fig:IEEE_37}Single-phase version of the IEEE-37 node distribution network from~\cite{baker2018network}. The orange and cyan dots represent the different locations of solar PV.}
\end{figure}

\begin{figure}[t]
    \vspace{-9pt}
  \subfloat [\label{fig:Pcomp_case1}]{   \includegraphics[width=0.48\linewidth,trim={0 0 0 0.5cm},clip]{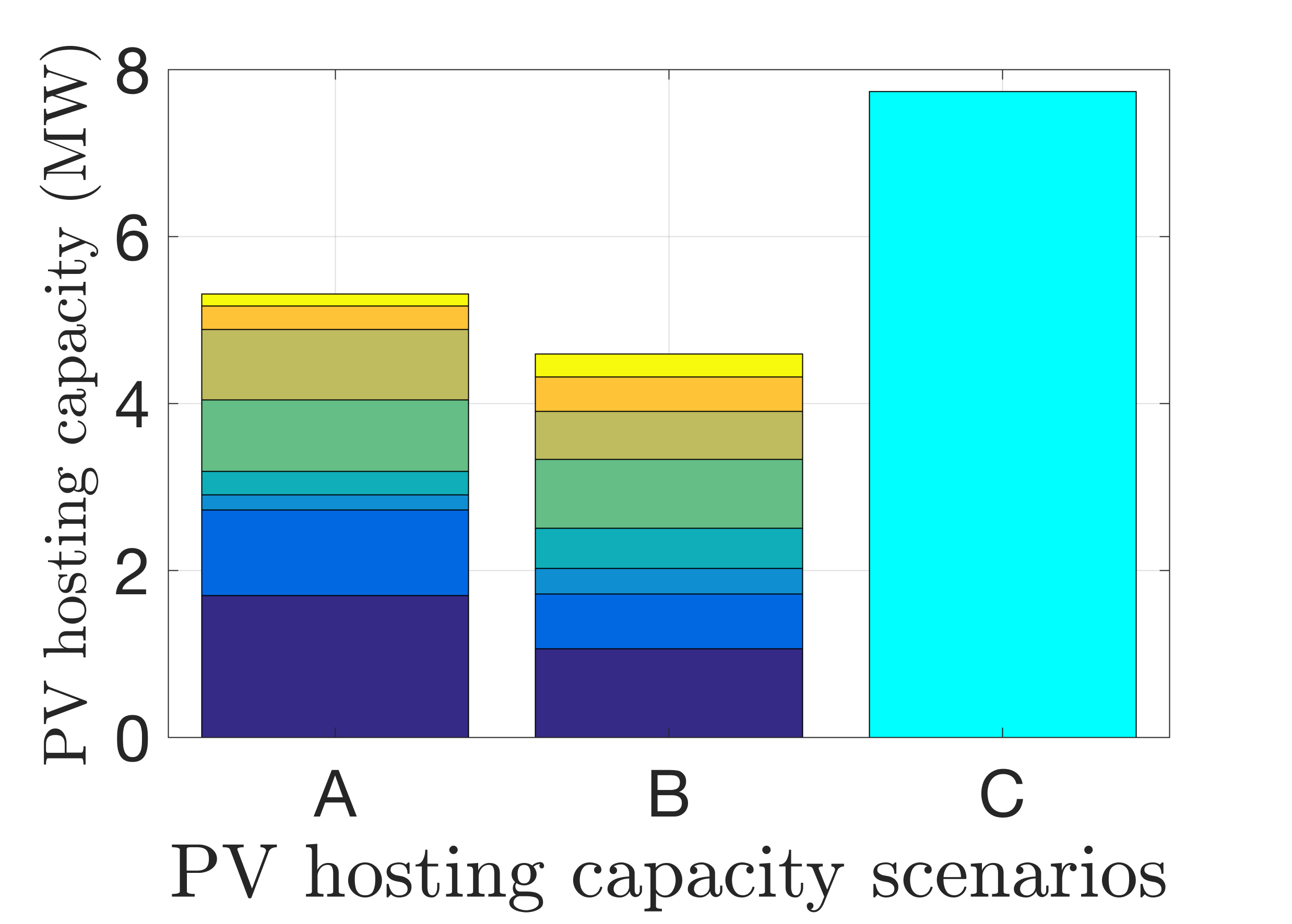}.}
    \hfill
  \subfloat [\label{fig:Vcomp_case1}]{    \includegraphics[width=0.48\linewidth,trim={0 0 0 0.5cm},clip]{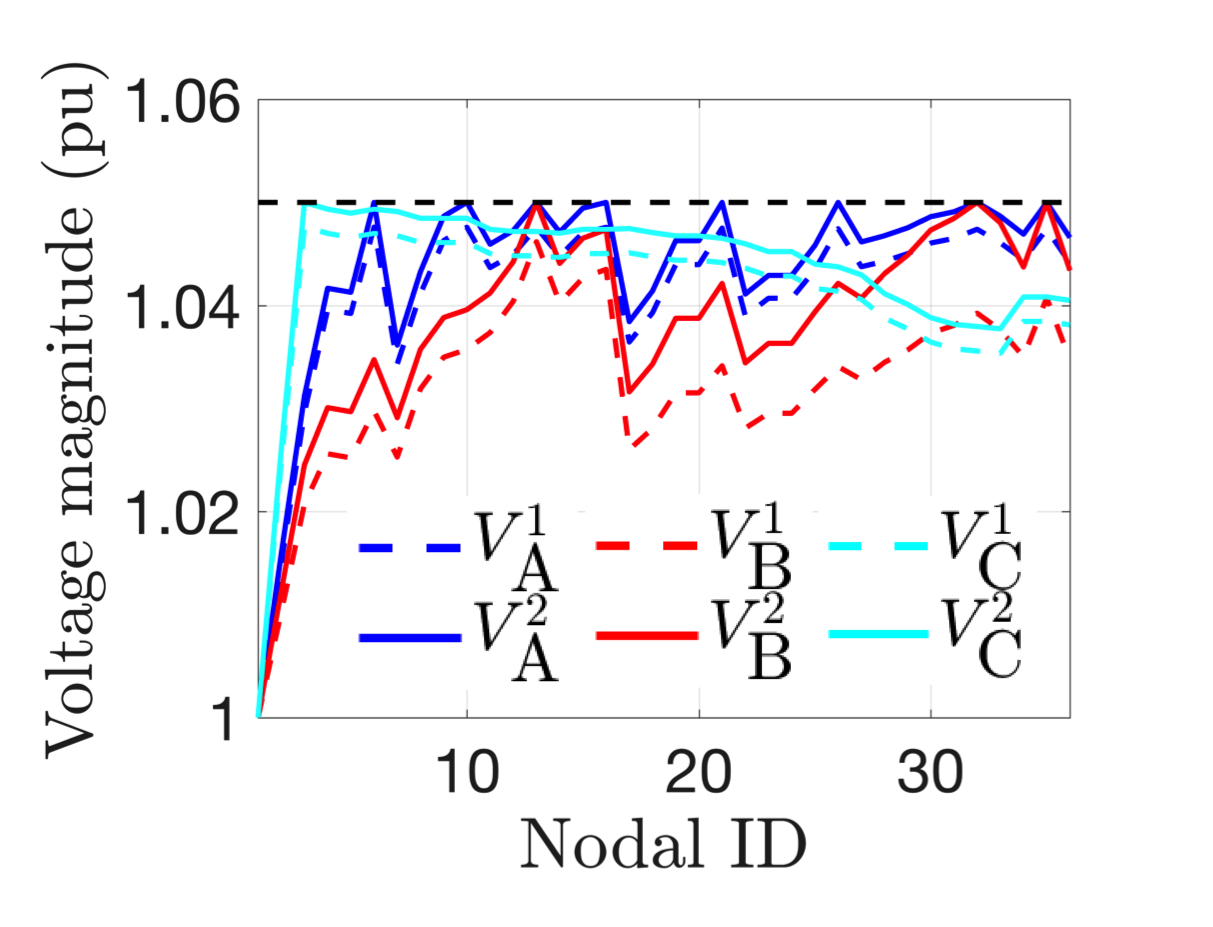}}
  \caption{Case study 1 on IEEE-37 node network for the three PV scenarios. (a) Shows the feeder's solar PV hosting capacity with Algorithm~\ref{alg1}.  (b) illustrates admissibility of PV hosting capacity via voltage profiles, where $V_{\text{A}}^1,V_{\text{B}}^1,V_{\text{C}}^1$ results from first iteration and $V_{\text{A}}^2,V_{\text{B}}^2,V_{\text{C}}^2$ are from the final iteration of Algorithm~\ref{alg1}.}
  \label{fig:case_study1}
\end{figure}

 Algorithm~\ref{alg1} is applied to  the IEEE-37 node distribution feeder shown in Fig.~\ref{fig:IEEE_37} for three different scenarios to determine $p_{\text{g},i}^+$. In this context, $p_{\text{g},i}^+$ can effectively be considered the solar PV hosting capacity.
The three different scenarios are specified in Table~\ref{table_CS1}. In scenarios~A (linear objective) and B (logarithmic objective), the solar PV units are installed at the (orange) leaf nodes with the largest demand, whereas in scenario~C (linear objective), solar PV is only allowed at (cyan) node~2 (e.g., utility-scale solar PV array). The optimization problem (P1) is solved with Gurobi~9.1 in Julia~1.1 in less than 1 sec and the solution is validated with Matpower~\cite{zimmerman2011matpower} on a standard MacBook Pro laptop with 2.2GHz CPU and 16GB RAM. The comparison of the resulting solar PV hosting capacity from each scenario using Algorithm~\ref{alg1} is shown in~Fig.~\ref{fig:Pcomp_case1}, with the stacked bars showing the allocated hosting capacity at the different nodes in the system. 
It can be seen that having a single centralized solar unit allows greater total solar PV capacity as compared to the distributed cases. The reason for this is that Scenario~C has fewer network limits to consider than the distributed case. This can be seen in Fig.~\ref{fig:Vcomp_case1}, which shows nodal voltages resulting from PV injections at the hosting capacity values after the first and second iteration of Algorithm~\ref{alg1}. As can be seen from the figure, voltages are at their upper limit at multiple nodes for scenarios A and B, but only at the head-node (node $2$) for scenario C. As a result, the distributed case (case A and B) has more active constraints and, hence, represents a more conservative solution as compared to the case with a single central PV (scenario C). Furthermore, scenario~B favors a more equitable allocation (log objective) that results in smaller net solar PV capacity ($\sum_i p_{\text{g},i}$), leading to reduced overall hosting capacity as compared to scenario~A.
Thus, Case Study~1 illustrates AC admissibility of the PV hosting capacity solutions resulting from Algorithm~\ref{alg1} and the effects of different objective function terms on the allocation of hosting capacity. The objective function design represents a powerful, but simple approach to align  hosting capacity with available incentives that prioritize certain load pockets or to ensure an equitable distribution of dispatchable DERs.

\begin{table}[t] \label{table_HCscenarios}
\centering
\caption{\label{table_CS1}PV hosting capacity scenarios}
{
\begin{tabular}{rll}
\toprule
{Scenario} & {Nodes with PVs} & {Objective function} \\
\midrule
            A & $\{5,9,12,15,20,25,31,34\}$ & $f_i(p_{\text{g},i}) = -p_{\text{g},i}$\\
            B & $\{5,9,12,15,20,25,31,34\}$ & $f_i(p_{\text{g},i}) = -\log(p_{\text{g},i})$\\
            C & $\{2\}$ & $f_i(p_{\text{g},i}) = -p_{\text{g},i}$
\\ \bottomrule
\end{tabular}
}
\end{table}

\subsection{Quantifying conservativeness}\label{sec:conservative}
In this section, we present simulation-based analysis to quantify the conservativeness of the convex inner approximation. {\color{black} We show how Algorithm~\ref{alg1} iteratively enlarges $\Delta p_{\text{g}}$ for the IEEE-13 node and IEEE-37 node test feeders. The results are tabulated in Table~\ref{table_iter}  and quantify how Algorithm~\ref{alg1} increases nodal capacity uniformly across the three scenarios (A, B, and C) and reduces the conservativeness of the CIA-based method.}

	\begin{table}[h!]
\centering
\caption{Feeder solar PV hosting capacity increase with iterations in Algorithm~1 for IEEE-13 node and IEEE-37 node systems.}\label{table_iter}
\begin{tabular}{cc}
IEEE-13 node system & IEEE-37 node system \\
\begin{tabular}{rlll}
\toprule
{Iteration} & {A} & {B} & {C} \\
\midrule
            Iter 1 (MW) & 6.6 & 6.3 & 8.8 \\
            Iter 2 (MW) & 9.1 & 8.6 & 12.1 \\
            \bottomrule
\end{tabular} &
\begin{tabular}{rlll}
\toprule
{Iteration} & {A} & {B} & {C} \\
\midrule
            Iter 1 (MW) & 5.1 & 4.4 & 7.4 \\
            Iter 2 (MW) & 5.3 & 4.6 & 7.7 \\
            \bottomrule
\end{tabular}\\
\end{tabular}
\end{table}

{\color{black} {Furthermore, to understand conservativeness of the CIA-based method, we compare it with other techniques from literature, including an SOCP relaxation and a non-convex, non-linear program (NLP). For the convex relaxation, we consider the SOCP formulation of the AC power flow equations and solve with GUROBI. The NLP formulation is based on the original non-convex AC OPF formulation and is solved using IPOPT. In general, the SOCP relaxation represents an upper bound on the nodal capacity, but can provide solutions that are not physically realizable, i.e., will violate network constraints. The non-convex NLP's (locally optimal) solutions, $p_{\text{g}}^-$ and $p_{\text{g}}^+$, are only guaranteed to be AC admissible at $p_{\text{g}}^-$ and $p_{\text{g}}^+$ and do not provide guarantees within the corresponding hyper-rectangle, $\Delta p_{\text{g}}$. The proposed CIA-based approach in (P1), however, guarantees not only that voltages and currents are within their limits over the entire range, but also that the nodal injections can be manipulated independently. The results in Table~\ref{table_comp} provide a measure of the conservativeness of the convex inner approximation by quantifying $p_{\text{g}}^-$ and $p_{\text{g}}^+$ for three distribution systems (13 node, 37 node and 123 node). The results indicate that the CIA-based solution is close to the NLP solution and, thus, is not overly conservative. The SOCP results provide much larger $p_{\text{g}}^-$ and $p_{\text{g}}^+$, however, they are not physically realizable, i.e., they violate network constraints. This comparison shows the effectiveness of the proposed approach in obtaining nodal capacity limits while guaranteeing satisfaction of network constraints.}
\textcolor{black}{Note that the CIA-based method takes less than ten seconds to determine $\Delta p_{\text{g}}$ for the IEEE-123 node system on a standard MacBook Pro with 2.2GHz CPU and 16GB RAM.}
\begin{table}[h!]
\centering
\caption{Comparing nodal capacities with the convex inner approximation (CIA), convex relaxation (CR), and non-convex (NLP) formulations}\label{table_comp}
{
\begin{tabular}{rlll}
\toprule
{System} & {CIA (MW)} & {NLP (MW)} & {CR (MW)} \\
\midrule
            13-node & [-1.5, 9.1] & [-1.5, 9.7] & [-1.5, 12] \\
            37-node & [-2.7, 5.3] & [-2.7, 5.3] & [-2.7, 16] \\
            123-node & [-4.5, 13.9] & [-4.5, 14] & [-4.5, 24]\\
            \bottomrule
\end{tabular}
}
\end{table}
}

\begin{remark}[Adapting analysis to distribution planning]
It is important to note that the nodal capacity herein can incorporate both flexible supply ($p_{\text{g}}^+>0$) and demand ($p_{\text{g}}^-<0$), but that $\Delta p_{\text{g}}$ is with respect to a particular operating point, $(P_{\text{L}},Q_{\text{L}})$. This is different from conventional PV hosting capacity studies that consider a representative annual, hourly demand profile~\cite{HOROWITZ2020106222}. In future work, we will extend (P1) and Algorithm~\ref{alg1} for multi-hour planning problems and incorporate battery storage and flexible demand to determine the ``dynamic hosting capacity'' of a feeder from quasi-static time-series (QSTS) analysis.
\end{remark}

{\color{black}\subsection{Existence and uniqueness of power flow solution}\label{sec:existence}
If a solution to the power flow equations in~\eqref{eq:dist_flow} exists, the  convex inner approximation (CIA) approach  presented herein guarantees satisfaction of network constraints for all, $p_{\text{g}}\in\Delta p_{\text{g}}$. However, the CIA-based approach in (P1) does not guarantee \textit{existence} of a power flow solution. Early work on the existence of power flow solutions in radial distribution networks showed that a unique solution exists for a wide range of practical parameter values~\cite{chiang1990existence}. However, there are power injections for which a solution may not exist, but such cases often occur under impractical operating conditions~\cite{molzahn2014investigation}. Nonetheless, ensuring existence of solutions for the entire range of nodal capacities is valuable. Recent works in literature, e.g.,~\cite{existance_PF1}, have provided sufficient conditions for the existence and uniqueness of a solution to~\eqref{eq:dist_flow}. This has since been extended to multi-phase distribution networks~\cite{bernstein2018load}. The existence of power flow solutions is also closely related to the voltage collapse problem and has been studied in~\cite{cui2019solvability}.  In this section, we utilize sufficient conditions from~\cite{existance_PF1}, to provide guarantees for both the satisfaction of network constraints and the existence of a solution over the range of nodal capacities. In this regard, we augment Algorithm~1 with two additional conditions that are adapted from~\cite{existance_PF1}. The first condition ensures a solution exists at the current operating point of Algorithm~\ref{alg1}, i.e., $x^0_{ij}=\text{col}\{P_{ij}^0,Q_{ij}^0,v_j^0\} \enspace \forall (i,j)\in \mathcal{L}$. This condition is given by:
\begin{align}
    \zeta(\hat{s})<u_{\min}^2\label{eq:EX_1}
\end{align}
where $\zeta(\hat{s})$ and $u_{\min}$ are defined as:
\begin{align}
    \zeta(s):=&||W^{-1}Y_{\text{LL}}^{-1}\overline{W}^{-1}\text{diag}(\overline{s})||_{\infty}\\
u_{\min}:=&\min_j|V_j^0/w_j|
\end{align}
	with $W:=\text{diag}(w)$, $w:=-Y_{\text{LL}}^{-1}Y_{\text{L}0}$,\, $s=p+\mathbf{i}q$ is the complex nodal power injection, $\hat{s}=-P_{\text{L}}-\mathbf{i}Q_{\text{L}}$, and $Y$ is the admittance matrix such that
		\begin{align*}
		    Y=\begin{bmatrix}Y_{00} & Y_{0\text{L}}\\
		    Y_{\text{L}0} & Y_{\text{LL}}\end{bmatrix}.
		\end{align*}
	Condition~\eqref{eq:EX_1} can be readily checked in Algorithm~\ref{alg1} at the operating point before each iteration. 
		
	 The second condition is incorporated into (P1) as follows. Define $s=\hat{s}+s_{\text{g}}$, where $s_{\text{g}} = p_{\text{g}} + \mathbf{i}q_{\text{g}}$. We can then determine sufficient conditions for existence of a solution over the range of nodal capacities $p_{\text{g}}^-\le p_{\text{g}}\le p_{\text{g}}^+$ with the following constraint added to (P1):
		\begin{align}
		    \Delta &:=\left(u_{\text{min}}-\frac{\zeta(\hat{s})}{u_{\min}}\right)^2-4 \zeta(s_{\text{g}})>0. \label{eq:EX_2}
		\end{align}
		To incorporate~\eqref{eq:EX_2} effectively, we define $\chi:=(u_{\text{min}}-\frac{\zeta(\hat{s})}{u_{\min}})^2/4$ and $A^{\text{w}}:=W^{-1}Y_{\text{LL}}^{-1}\overline{W}^{-1}$, which yields condition $||A^{\text{w}}\text{diag}(s_{\text{g}})||_{\infty}<\chi$,
	from which we employ the definition of matrix norm to get expression
		\begin{align}
		 \max_{i=1,\hdots,N} \sum_{j=1}^N|A^{\text{w}}_{ij}\overline{s}_{\text{g},j}|<\chi . \label{eq:EX_5}
		\end{align}
			The constraint in~\eqref{eq:EX_5} can also be represented by $N$ constraints of the form $\sum_{j=1}^N|A^{\text{w}}_{ij}\overline{s}_{\text{g},j}|< \chi \enspace \forall i=1,\hdots, N$.
		By defining $A^{\text{w}}_{ij} = a^{\text{w}}_{ij} + \mathbf{i} b_{ij}$ and using and expanding the complex product, we get the following equivalent convex formulation of~\eqref{eq:EX_2} that is composed of $N$ linear inequalities and $N$ second-order cone (SOC) constraints:
\begin{align}
\label{eq:EX_7c3} \tag{C3}
  \begin{gathered}
\sum_{j=1}^N t_{ij} < \chi \quad \forall i=1,\hdots, N  \\ 
\left|\left|
\left[
\begin{matrix}
a^{\text{w}}_{ij} & b^{\text{w}}_{ij} \\
b^{\text{w}}_{ij} & -a^{\text{w}}_{ij}
\end{matrix}
\right]
\left[ \begin{matrix}
p_{\text{g},j} \\
q_{\text{g},j}
\end{matrix}
\right]
\right|\right|_2 \le t_{ij}\quad \forall j=1,\hdots, N .
\end{gathered}
\end{align}

		The convex reformulation in~\eqref{eq:EX_7c3} can now be readily included in (P1) to compute the nodal capacities with guarantees on the existence of an AC power flow solution. To investigate the conservativeness of the CIA-based method with and without~\eqref{eq:EX_7c3},  Table~\ref{table_exis} compares the hosting capacity  for the 13-, 37-, and 123-node networks. The results indicate that including~\eqref{eq:EX_7c3} still leads to meaningful and practical solutions similar to the results in Section~\ref{sec:conservative}.
		
		\begin{table}[h!]
\centering
\caption{Comparing proposed convex inner approximation with and without existence condition (C3)}\label{table_exis}
{
\begin{tabular}{rlll}
\toprule
{Type} & {13-node} & {37-node} & {123-node} \\
\midrule
            Without C3 (MW) & [-1.5, 9.1] & [-2.7, 5.3] & [-4.5, 13.9] \\
            With C3 (MW) & [-1.5, 8.8] & [-2.7, 5.3] & [-4.5, 13.8] \\
            \bottomrule
\end{tabular}
}
\end{table}

		}
The analysis and simulation results presented in Case Study~1 and this section set $q_{\text{g}}=0$, e.g., no controllable reactive power. However, the role of reactive power management in optimizing DER nodal capacities is valuable and, thus, the focus of the next section.

\section{Role of reactive power}\label{sec:Reac_power} 
Controlling reactive power, $q_{\text{g}}$, can increase the nodal capacity, $\Delta p_{\text{g}}$ as shown by analyzing the effect of different reactive power control schemes in this section. Specifically, we will compare between DERs that are operated at unity power factor, fixed power factor, and those with advanced reactive power control capability, such as IEEE Standard 1547~\cite{1547}. The different reactive power schemes along with the relevant relations between $q_\text{g}$ and $p_\text{g}$ are provided in Table~\ref{table_scenarios}. For each particular scheme, the corresponding constraints are added to (P1) when determining the nodal capacity.

\begin{figure}[t]
    \vspace{-9pt}
  \subfloat [\label{fig:reac_cap_comp}]{   \includegraphics[width=0.48\linewidth,trim={0 0 0 0.5cm},clip]{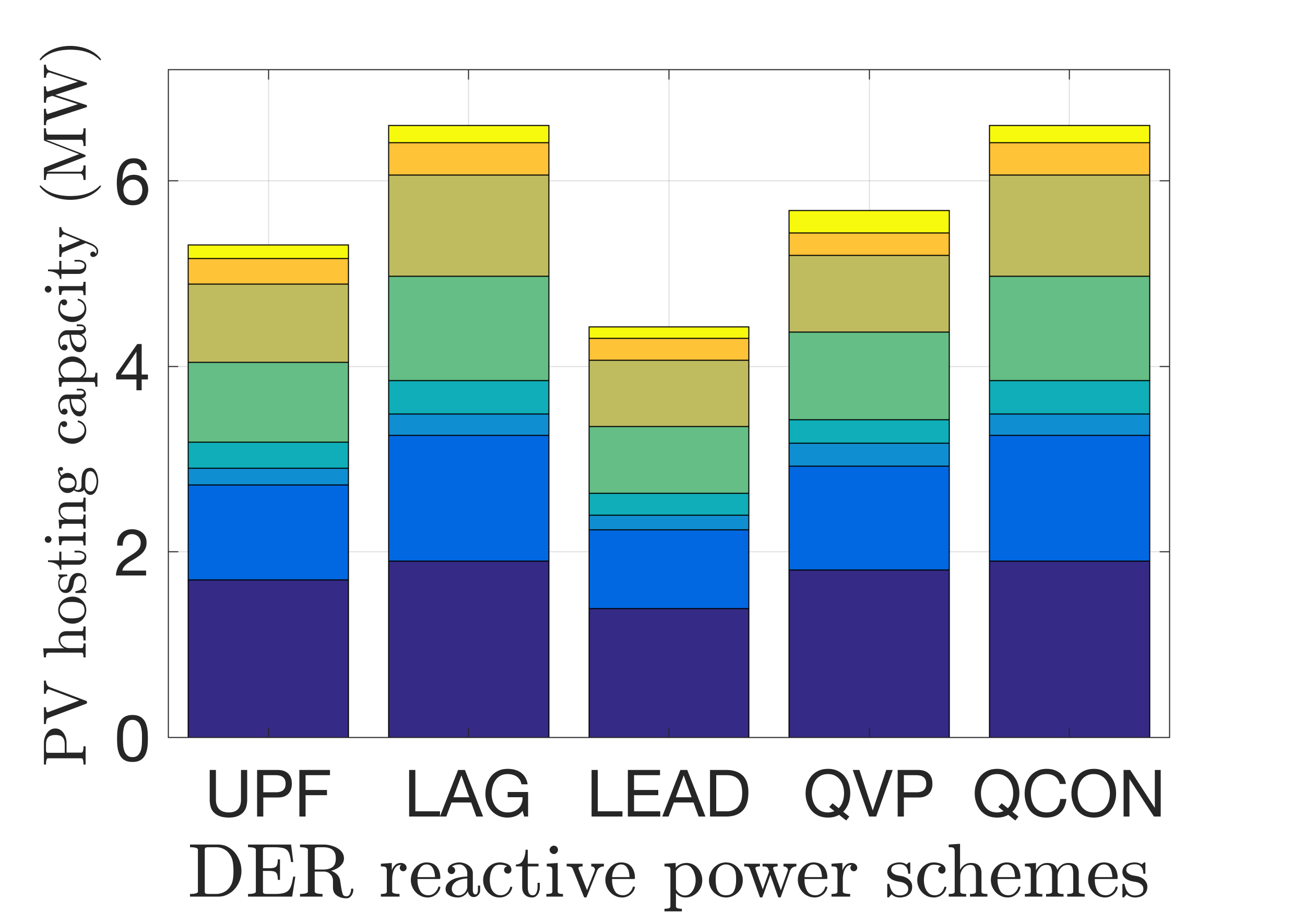}}
    \hfill
  \subfloat [\label{fig:reac_volt_comp}]{    \includegraphics[width=0.48\linewidth,trim={0 0 0 0.5cm},clip]{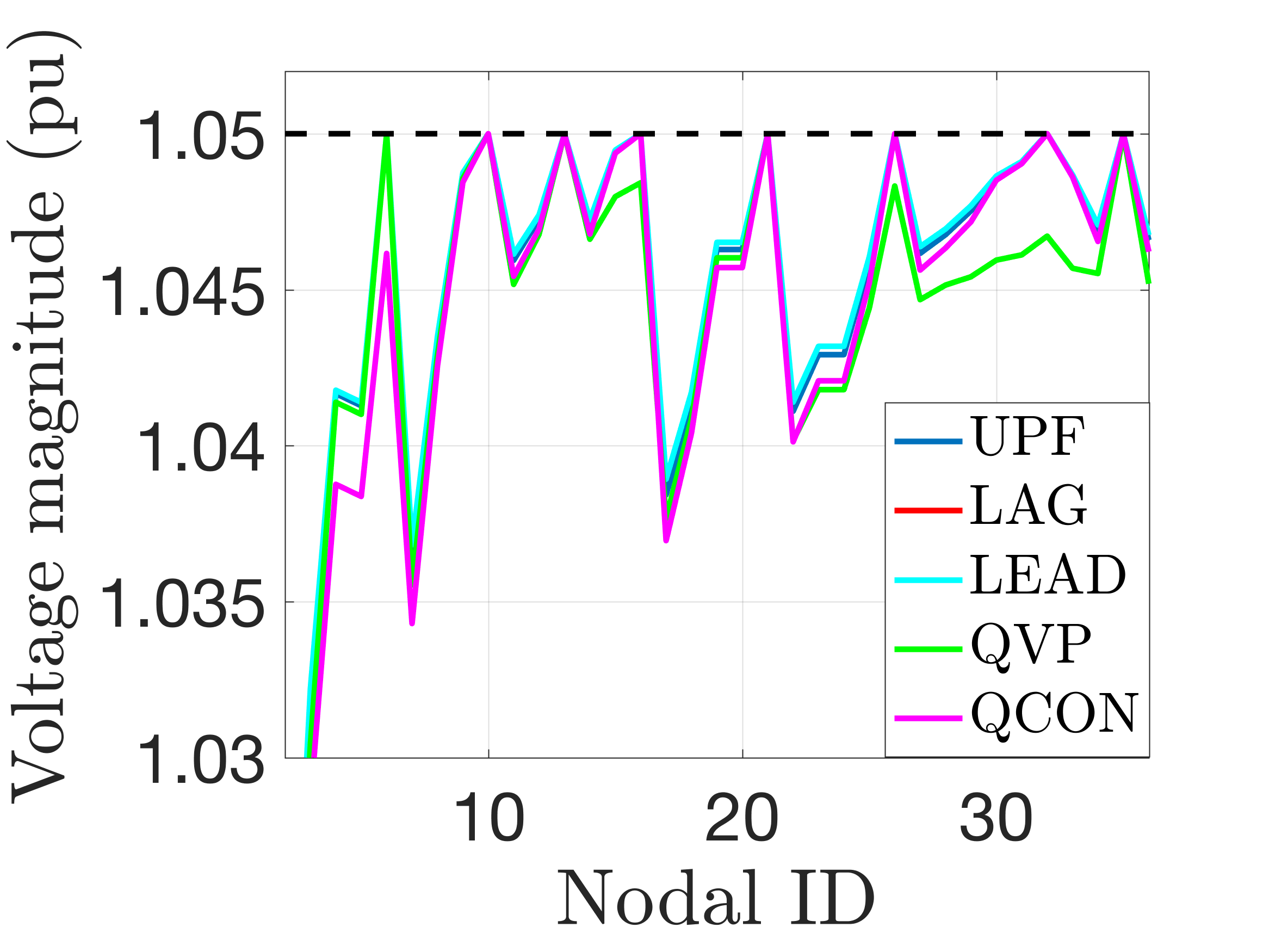}}
  \caption{Reactive power schemes for case Study~1 (Scenario~A) on the IEEE-37 node network for five different reactive power schemes: (a) Solar PV hosting capacity for each reactive power scheme after employing Algorithm~\ref{alg1} (b) Illustrating admissibility with voltage profiles for the final iterate from Algorithm~\ref{alg1}.}
  \label{fig:case_study2}
\end{figure}

Fig.~\ref{fig:reac_cap_comp} compares the feeder's solar PV hosting capacities, $\sum_i p_{\text{g},i}^+$, resulting from the different reactive power schemes applied to Scenario~A of Case Study~1. The stacked bar chart in Fig.~\ref{fig:reac_cap_comp} also shows the hosting capacity at the different nodes with DERs in this system. Scheme~UPF represents the hosting capacity with unity power factor, which matches the result from Scenario~A in Fig.~\ref{fig:Pcomp_case1} and serves as the base-case for comparison. Scheme~LAG employs a lagging power factor of 0.95 ($\gamma_i=-0.33$), while LEAD uses a leading power factor of 0.95 ($\gamma_i=+0.33$). Scheme QVP employs a common volt-VAr policy with $\beta_i^0=0$ and $\beta_i^1=-0.073$, while QCON represents advanced inverter capability with quadratic constraints and $\bar S_{\text{g},i}=2$MVA and a minimum power factor of 0.95. The results show that for scheme LEAD, the hosting capacity is reduced while schemes LAG, QCON and QVP increase hosting capacity. In LEAD, this is due to reactive power injections increasing with active power injections resulting in larger $v$ and, hence, reduces $p_{\text{g}}^+$. The opposite occurs in the other schemes. Interestingly, QCON achieves the same nodal capacity as LAG by bringing power factors to their lower limit of 0.95 in order to maximize nodal capacity. The voltage profiles at the hosting capacities for the different schemes are compared in Fig.~\ref{fig:reac_volt_comp} and are clearly AC admissible. 
\textcolor{black}{ The results in Fig.~\ref{fig:case_study2} highlight the complex relationship between PV hosting capacity and reactive power control and how CIA-based methods can be used to effectively study siting, sizing, and dispatch of DERs in a grid-aware manner. 
This is in line with the ongoing developments of reactive power requirements and standards~\cite{1547}.
}

\begin{table}[h!] \label{table_Qschemes}
\centering
\caption{\label{table_scenarios}DER reactive power schemes}
{
\begin{tabular}{rll}
\toprule
{Scheme} & {Description} & {Constraint $(g_i(p_{\text{g},i},q_{\text{g},i},v_i))$} \\
\midrule
            UPF & Unity power factor & $q_{\text{g},i}=0$\\
            LAG & Lagging power factor & $q_{\text{g},i}=-\gamma_i p_{\text{g},i}$\\
            LEAD & Leading power factor & $q_{\text{g},i}=\gamma_i p_{\text{g},i}$\\
             QVP & Volt-VAr policy & $q_{\text{g},i}=\beta_i^0+\beta_i^1v_i$\\
            QCON & Quadratic constraint & $p_{\text{g},i}^2+q_{\text{g},i}^2\le \overline{S}_{\text{g},i}^2$
           
\\ \bottomrule
\end{tabular}
}
\end{table}

The next section employs the nodal capacity, $\Delta p_{\text{g}}$, 
to develop a simple, open-loop, decentralized DER control policy for the realtime, grid-aware disaggregation of a power reference signal. This turns the whole feeder into a responsive grid resource with \textit{a-priori} AC admissibility guarantees that can provide fast grid services. 

\section{Realtime grid-aware  disaggregation}\label{sec:real_time}
Dispatching a set of networked DERs in response to a fast, time-varying power reference signal while guaranteeing admissible operations is a challenging problem. 
However, it is necessary to solve this problem before aggregators can safely coordinate millions of behind-the-meter DERs without jeopardizing reliability of the grid.
Thus, after computing the available nodal capacity (offline), as shown in Fig.~\ref{fig:cyber_physical_model}, this section proposes a simple, grid-aware controller  to allocate the required flexibility among the available resources in the network (i.e., disaggregate the signal) in realtime. 
\textcolor{black}{ The advantage of such a mechanism is that it simplifies the interface between the distribution system operator (DSO) and the individual aggregators. The DSO has the responsibility of determining the network's DER nodal capacities, whereas the aggregator is only required to operate within the DSO's provided nodal capacities. This is a reasonable setting as the DSO has access to the network data while the aggregators may not, and aggregators have access to DERs while DSO's may not. Hence, the proposed CIA-based method enables a bridge between the DSO and aggregator, taking advantage of each of their capabilities and roles, in order to enable grid-aware, real-time dispatch of DERs.}


The necessary parameters to execute the realtime, grid-aware disaggregation are $p_{\text{g},i}^+$ and $p_{\text{g},i}^-$ and can be updated every 15-60 minutes by the grid operator executing Algorithm~\ref{alg1}, which is the timescale of the baseline of the aggregate uncontrollable net-demand.
The realtime disaggregation can then be solved by a DER aggregator to provide fast grid services without the need to include \textit{any} information about the underlying grid parameters. That is, the nodal capacities embed the AC OPF parameters and constraints to simplify the aggregator's dispatch. {\color{black} The disaggregation process is determined by the following open-loop policy executed by the aggregator at each node $i$ and discrete-time step $k$:
\begin{align}
    p_{\text{g},i}[k]=
    \begin{cases}
    \min\{\frac{p_{\text{g},i}^+}{\sum_ip_{\text{g},i}^+}P_{\text{ref}}[k],p_{\text{g},i}^+\} \qquad P_{\text{ref}}[k]\ge0\\
   \max\{\frac{p_{\text{g},i}^-}{\sum_ip_{\text{g},i}^-}P_{\text{ref}}[k],p_{\text{g},i}^-\} \qquad P_{\text{ref}}[k]<0
    \end{cases}
    .
    \label{eq:disagg_RT}
\end{align} }  
The next case study shows the effectiveness of the proposed disaggregation process in having DERs collectively respond to grid service signals while guaranteeing AC admissibility. {\color{black}Future work will consider the role of feedback and disturbances.}

\subsubsection*{\textbf{Case study 2}} 
\begin{figure}[t]
    \vspace{-9pt}
  \subfloat [\label{fig:head_tracking}]{   \includegraphics[width=0.48\linewidth,trim={0 0 0 0.5cm},clip]{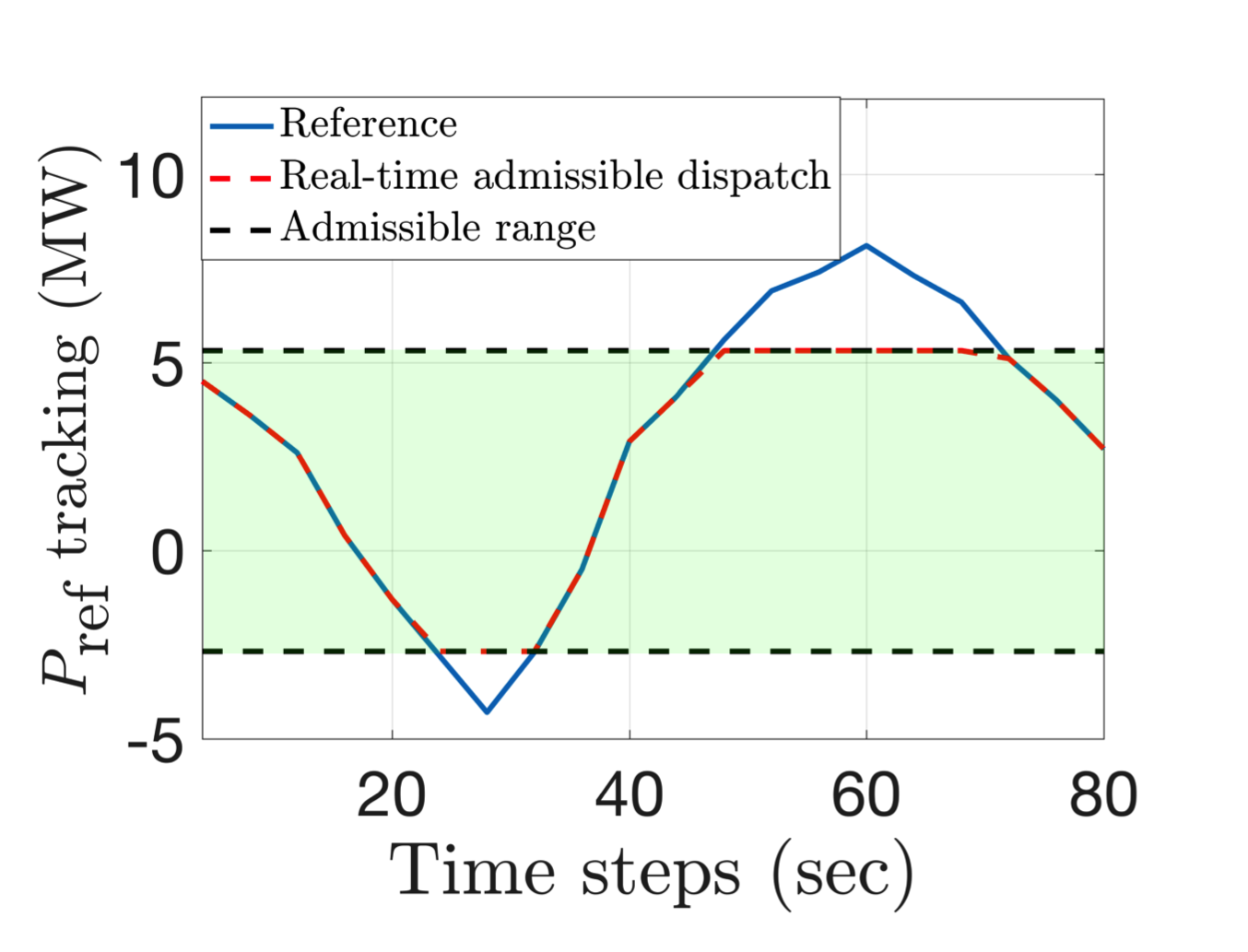}}
    \hfill
  \subfloat [\label{fig:V_limit}]{    \includegraphics[width=0.48\linewidth,trim={0 0 0 0.5cm},clip]{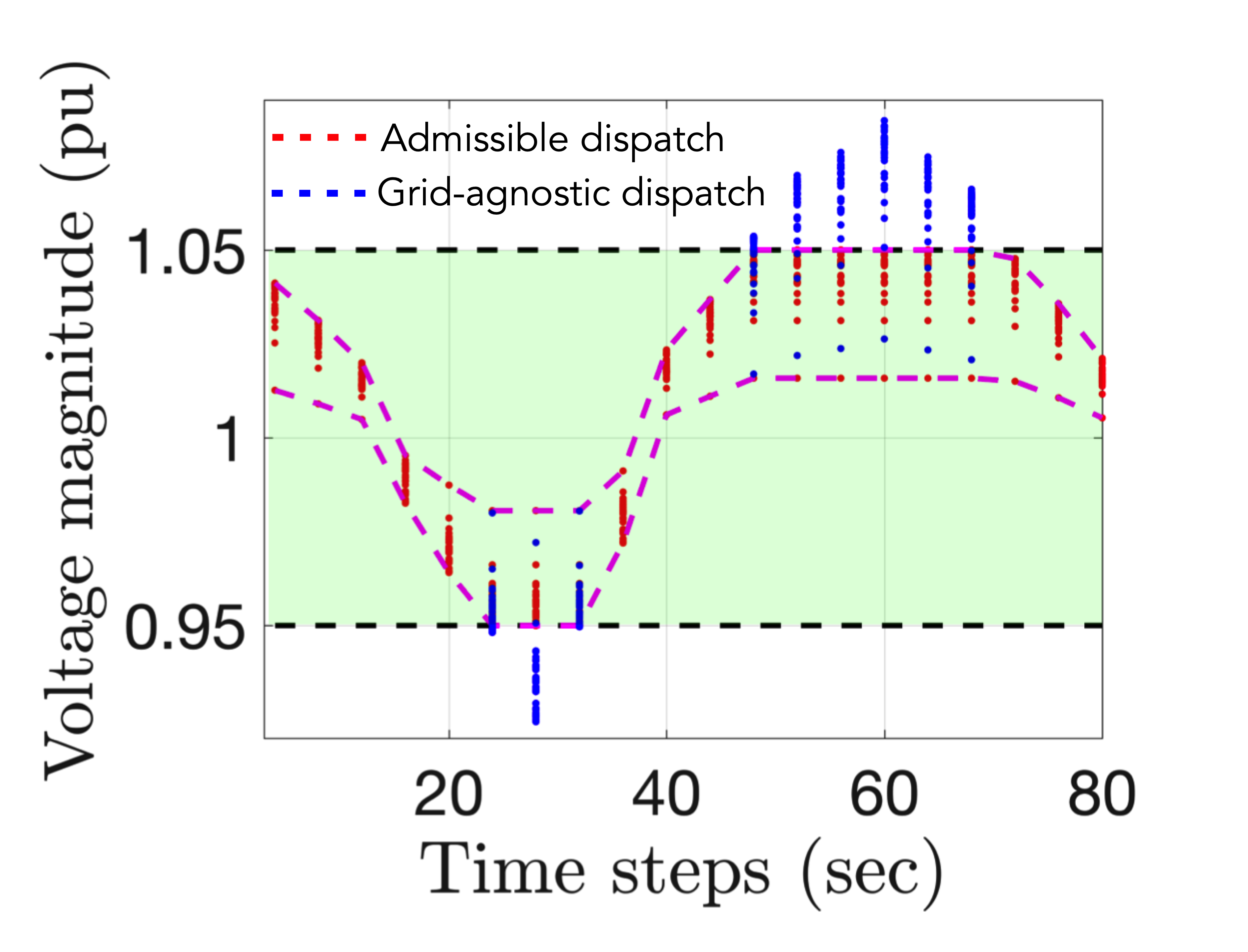}}
  \caption{Case study 2: (a) tracking performance of the realtime disaggregation policy shown in~\eqref{eq:disagg_RT} for IEEE-37 node system (b) Voltage profile of the IEEE-37 node network showing admissible dispatch when following the grid-aware disaggregation policy ({\color{red}red}) and voltage violations when following a grid-agnostic approach ({\color{blue}blue}). Note that the envelope of admissible voltage magnitudes from the grid-aware dispatch is shown in \textcolor{mypink1}{pink}. The grid-agnostic approach results in a maximum voltage violation of $0.03$pu at time-step $60$s.}
  \label{fig:tracking}
\end{figure}

The effectiveness of the offline Algorithm~\ref{alg1} and the online disaggregation in~\eqref{eq:disagg_RT}
is illustrated with the IEEE-37 node system, where we use the nodal capacities defined by Scenario~A. The case study shows that the feeder is being managed within its limits at all times despite providing a large range of flexibility from the responsive DERs. Fig.~\ref{fig:head_tracking} shows a power reference grid service signal and the aggregate response from dispatching the DERs. As shown, the reference signal is tracked well when it is within the admissible range and the aggregator dispatch is also AC-admissible as shown in Fig.~\ref{fig:V_limit}. In a practical setting, the DER aggregator should only offer what can be delivered, but the case study is meant to illustrate $i$) how the disaggregation enables realtime, grid-aware nodal dispatch and $ii$) how the nodal capacity, $\Delta p_{\text{g}}$, enables simple aggregation to characterize a radial feeder's admissible range. 
  Clearly, if the aggregator was grid-agnostic and just coordinated DERs to ensure perfect tracking,  
  then 
  such a ``greedy'' version of the realtime DER dispatch 
  leads to network voltages violations, as illustrated by the blue dots in Fig.~\ref{fig:V_limit}. Thus, the proposed open-loop control scheme is grid-aware and scalable across a network of DERs by just broadcasting a single scalar power reference signal.

\section{Conclusions and Future work}\label{sec:conclusion}
This manuscript presents 
a convex inner approximation of the AC OPF problem. Leveraging convex lower and upper bounds on the nonlinear branch flow terms in the AC formulation, the inner approximation enables a characterization of nodal capacity that is provably AC admissible across the entire range. A novel algorithm is then presented to successively improve the nodal capacity of a feeder. The effect of different reactive power control schemes is then illustrated and volt-VAr and advanced inverter schemes are shown to further improve the nodal capacity. Finally, a realtime disaggregation scheme is developed for dispatching flexible demand in realtime across the network, while respecting the grid constraints, thus, enabling grid-aware fast grid services.

Future work will extend this work to multi-phase feeder models and meshed networks to account for more realistic distribution feeders. Optimizing legacy and other front-of-meter grid asset schedules to increase or maintain the nodal hosting capacities is also of interest. Finally, employing feedback from salient grid measurements to provide robust admissibility guarantees in realtime under changes to expected demand/solar PV is another important area to investigate. 

\bibliographystyle{IEEEtran}
\IEEEtriggeratref{35}
\small\bibliography{sample.bib}

\begin{IEEEbiography} [{\includegraphics[width=1in,height=1.25in,clip,keepaspectratio]{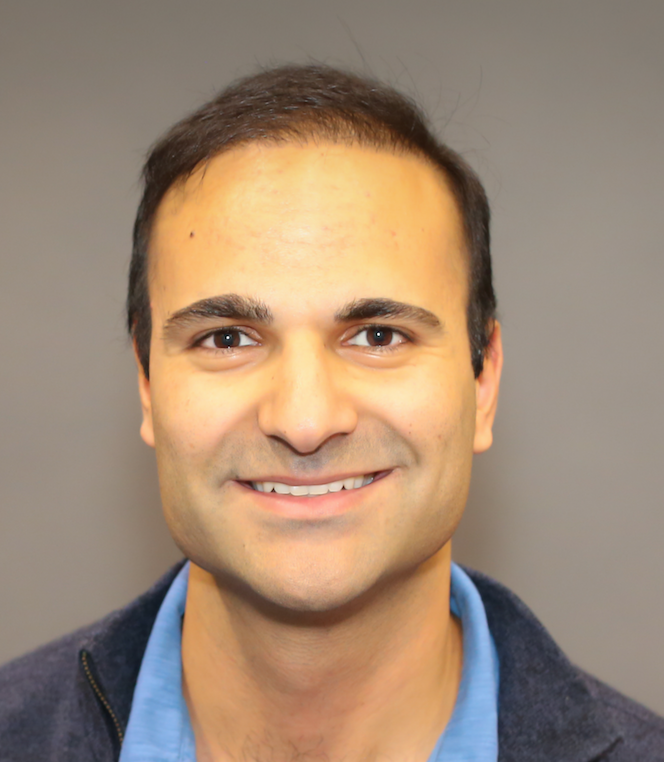}}] {Nawaf Nazir}(S’17-M'21) received the
M.S. degree in electrical engineering from Virginia Polytechnic Institute and State University, Blacksburg, VA, USA, in 2015 and the PhD degree in Electrical Engineering from the Department of Electrical and  Biomedical Engineering at the University of Vermont, Burlington, VT, USA. He is currently a Postdoctoral Researcher at the Pacific Northwest National Lab, Richland, WA, USA.
His research interests include optimization, control and machine learning applied to complex networked systems, and, in particular, emphasizes reliability, resilience and real-time control of energy systems.
\end{IEEEbiography}

\begin{IEEEbiography}[{\includegraphics[width=1in,height=1.25in,clip,keepaspectratio]{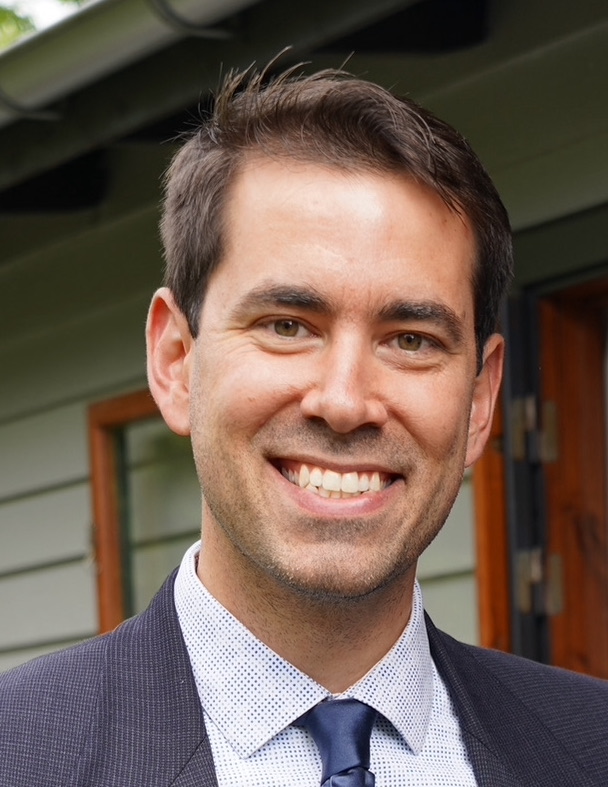}}]{Mads Almassalkhi}(S'06-M'14-SM'19) is Associate Professor in the Department of Electrical and Biomedical Engineering at the University of Vermont, Burlington, Vermont, co-founder of startup Packetized Energy, and Chief Scientist at Pacific Northwest National Laboratory (by joint appointment). His research interests lie at the intersection of power systems, mathematical optimization, and control systems and focus on developing scalable algorithms that improve responsiveness and resilience of energy and power systems. He is currently serving as Chair of the IEEE Power \& Energy Society's Smart Buildings, Loads, and Customer Systems (SBLC) technical subcommittee on Loads. Prior to joining the University of Vermont, he was lead systems engineer at energy startup company Root3 Technologies in Chicago, IL. Before that, he received his PhD from the University of Michigan in Electrical Engineering: Systems in 2013 and a dual major in Electrical Engineering and Applied Mathematics at the University of Cincinnati in 2008.
\end{IEEEbiography}

\end{document}